\documentclass[a4paper,12pt,fleqn]{article}
\usepackage{amsmath,amssymb,amsthm,graphicx,subfigure,float,caption,epstopdf,tabularx,color, bm,amsfonts,epic}
\usepackage[top=1.25in, bottom=1.2in, left=1.0in, right=0.9in]{geometry}
\usepackage{appendix}
\usepackage{multirow}
\usepackage{xcolor}
\usepackage{setspace}

\allowdisplaybreaks
\newtheorem{thm}{Theorem}[section]
\newtheorem{lem}{Lemma}[section]

\newtheorem{rmk}{Remark}[section]
\newtheorem*{prf}{Proof}
\numberwithin{equation}{section}

\usepackage{indentfirst}
\graphicspath{{Fig}}
\begin{document}

\title{ A structure-preserving algorithm for the fractional nonlinear
Schr\"{o}dinger equation based on the SAV approach}

\author{Yayun Fu, \quad Wenjun Cai, \quad Yushun Wang\footnote{Correspondence author. Email:wangyushun@njnu.edu.cn.}\\{\small Jiangsu Key Laboratory for NSLSCS,}\\{\small School of Mathematical Sciences,  Nanjing Normal University,}\\{\small  Nanjing 210023, China}\\}
\date{}
\maketitle

\begin{abstract}

The main objective of this paper is to present an efficient structure-preserving scheme, which is
 based on the idea of the scalar auxiliary variable approach, for solving the space fractional nonlinear Schr\"{o}dinger equation. First, we reformulate the equation as a Hamiltonian system, and obtain a new equivalent system via introducing a scalar variable. Then, we construct a semi-discrete energy-preserving scheme by using the Fourier pseudo-spectral method to discretize the equivalent system in space direction. After that, applying the Crank-Nicolson method on the temporal direction gives a linear implicit scheme in the fully-discrete version. As expected, the proposed scheme can preserve the energy exactly and more efficient in the sense that only decoupled equations with constant coefficients need to be solved at each time step. Finally, numerical experiments are provided to demonstrate the effectiveness and conservation of the scheme.\\[2ex]

\end{abstract}

 %%%%%%%%%%%%%%%%%%%%%%%%%%%%%%%%%%%%%%%%%% begin section 1 introduction %%%%%%%%%%%%%%%%%%%%%%%%%%%%%%%%%%%%%%%%%%%%
\section{Introduction}

Due to the memory and genetic property of fractional operators, fractional differential equations are more suitable to model various scientific and engineering problems with long-range temporal cumulative memory effects and spatial interactions than integer order, thus they are widely implemented in the fields of biomedical engineering \cite{p1,p2}, physics \cite{p3}, and hydrological applications \cite{p4}. The space fractional nonlinear Schr\"{o}dinger (NLS) equation, introduced in Refs. \cite{p5, p6} via extending the Feynman path integral to the L\'{e}vy integral, is fractional version of classical NLS equation and consider long-range interactions. In quantum mechanics, the fractional NLS equation is more accurate than the integer order in describing how the quantum state of a physical system changes in time \cite{p8,p24}. Nowadays, the studies for the equation have been considered and many significant achievements have been made \cite{p7, p8, p22}.

In this paper, we numerically consider the fractional NLS equation as follows
\begin{align}\label{NLS:eq:1.1}
\text{i}\frac{\partial u(\bm x,t)}{\partial t}-\gamma(-\Delta)^{\frac{\alpha}{2}}u(\bm x,t)+(V(\bm x)+\beta|u(\bm x,t)|^2)u(\bm x,t)=0,\ \ t \in (0,T],
\end{align}
with the initial condition
\begin{align}\label{NLS:eq:1.2}
u(\bm x,0)=u_0(\bm x),
\end{align}
where $1< \alpha \leq 2$, $\text{i}^2=-1$, $u(\bm x,t)$ is a complex-valued wave function of $\bm x\in \mathbb{ R}^d$ ($d$=1, 2) and $u_0(\bm x)$ is a given smooth function, $V(\bm x)$ is an arbitrary potential function,
%and $\Omega=(a,b)\times (c,d)$ is a bounded domain in $\mathbb{ R}^2$. The unknown $u(x,y,t)$ is a complex-valued wave function, $u0(x,y)$ is a given smooth function.
$\gamma> 0$ is a small semiclassical parameter that describes the group velocity dispersion, $\beta$ is a real constant with positive value for focusing (or attractive) nonlinearity. The fractional Laplacian $-(-\Delta)^{\frac{\alpha}{2}}$
is defined as a pseudo-differential operator with the symbol $|\bm \xi|^{\alpha}$ in the Fourier space \cite{p13}
\begin{align}\label{eq:1.4}
-(-\Delta)^{\frac{\alpha}{2}}u(\bm x)=\mathcal{F}^{-1}[-|\bm \xi|^{\alpha}\mathcal{F}[u]],
\end{align}
where $\mathcal{F}$ is the Fourier transform and $\mathcal{F}^{-1}$ denotes its inverse. If $\alpha=2$, the fractional NLS equation reduces to the standard NLS equation.

Following a similar argument in Ref. \cite{p7}, we can prove that system \eqref{NLS:eq:1.1} with the periodic boundary condition possesses the following energy and mass conservation laws
\begin{align}\label{NLS:eq:1.3}
E(t)=E(0),\ \ M(t)=M(0),
\end{align}
where the energy is defined as
\begin{align}\label{NLS:eq:1.4}
E(t):=\int_{\mathbb{ R}^d} \Big[-\frac{\gamma}{2}|(-\Delta)^{\frac{\alpha}{4}}u({\bm x},t)|^{2}+\frac{1}{2}V(\bm x)|u({\bm x},t)|^{2}+\frac{\beta}{4}|u({\bm x},t)|^{4}\Big]d{\bm x},
\end{align}
and the mass has the form
\begin{align}\label{NLS:eq:1.5}
 M(t):=\int_{\mathbb{ R}^d}|u({\bm x},t)|^{2}d{\bm x}.
\end{align}

The prior researches generally confirmed that the structure-preserving methods conserving one or more intrinsic properties of a given dynamical system, are more superior than the traditional methods in long time stability for numerical simulations \cite {p30, p31, p32, p33}. Nowadays, investigating the structure-preserving numerical schemes for fractional equations \cite{p40,p41,p42}, especially for the fractional NLS equation \cite{p18,p19,p21,p23,p25, p27, p29, p34}, has captured researchers' increasing attention. The conservation of energy is a crucial property for the fractional NLS equation. In recent years, constructing various energy-preserving schemes for the fractional NLS equation exerted a tremendous fascination on scholars. For instance, in Ref. \cite{p18,p19,p21}, Wang and Xiao presented some energy-preserving difference schemes for the problem. Duo et al. \cite{p23} proposed a Crank-Nicolson Fourier spectral method for solving the equation. An energy conservative finite element scheme was studied by Li et al. in Refs. \cite {p29}. For two-dimensional case, some energy-preserving schemes were also constructed \cite{p24, p35}. %What is more, in order to avoid heavy iterative costs at each time step, a linearly implicit energy-preserving scheme for the fractional NLS equation is very necessary.
However,
%there exist few works focusing on studying the linear implicit energy-preserving scheme, just the Refs.\cite{p18, p29} for the one-dimensional case.
%However, as far as we know, there exist few works focusing on studying the linear implicit energy-preserving scheme, just the Refs.\cite{p18, p29} for the one-dimensional case.
 %What is more, in order to avoid heavy iterative costs at each time step, there exist few works focusing on studying the linear implicit energy-preserving scheme, just the Refs.\cite{p18, p29} for the one-dimensional case.
%Moreover, there are some works have been developed the energy-preserving scheme
%For two-dimensional case, some energy-preserving schemes were also constructed \cite{p24, p35}. What is more, in order to avoid heavy iterative costs at each time step, there exist some works focusing on studying the linear implicit energy-preserving scheme for the one-dimensional fractional NLS equation \cite{p18, p29}.
%However, as far as we know,
most of the energy-preserving schemes for the equation are fully-implicit, therefore one needs to use iterations to solve a system of nonlinear algebraic equations at each time step. As far as we know, there exists few works focusing on studying the linear implicit energy-preserving scheme for the fractional NLS equation. For example, Wang \cite{p18} proposed a linear implicit difference scheme for the equation, which conserves both the mass and energy. Li et al. \cite{p29} studied a linear finite element conservative scheme for the problem, and discussed the convergence of the scheme.

%The only research achievements are Wang and Li  \cite{p18, p29} considering linear difference and finite element conservative schemes for the fractional NLS equation, which are efficient and the algebraic equations with variable coefficients need to be solved at each temporal step.

%one needs to use iterations to solve a system of algebraic equations at each time step, there exist few works focusing on studying the linear implicit energy-preserving scheme, just the Refs. \cite{p18, p29} for the one-dimensional case.
%only a few studies have focused on investigating linear energy-preserving scheme for the one-dimensional fractional NLS equation \cite{p18,p19,p29}.

 The scalar auxiliary variable (SAV) approach was built upon the invariant energy quadratization (IEQ) method \cite{p36,p37}, which inherited all advantages of IEQ approach but also overcame most of its shortcomings, and has been used to construct efficient and accurate numerical schemes for a large class of gradient flows \cite {p38, p39,p46}. To the best of our knowledge, there is no previous research using the SAV approach to construct energy-preserving scheme for fractional differential equations. As we all known, the Fourier pseudo-spectral method is an important method to solve differential equations \cite{p45}, and its main properties are high accuracy and much more efficient.
 In this paper, to overcome the difficulty of the fully-implicit scheme and reduce the the computational complexity,
we develop a linear implicit energy-preserving scheme for the fractional NLS equation based on the the SAV approach and the Fourier pseudo-spectral method. Compared with the schemes in
Refs. \cite{p18, p29}, the proposed scheme is easy to implement and much more efficient,
for the reason that one only needs to solve decoupled linear equations with constant
coefficients at each time step.

 %As we all known, the spectral methods are a class of high accuracy and much more efficient methods, and are used to solve differential equations \cite{p45}. Therefore the Fourier pseudo-spectral method is used to discrete the equivalent system in space direction to reduce the the computational complexity in the implementation by introducing the Fourier transform technique. Therefore, the Fourier pseudo-spectral method is used to discrete the equivalent system in space direction to reduce the the computational complexity in the implementation by introducing the Fourier transform technique. The spectral method is a fast, high order and widely used technique to solve differential equations \cite{p45}. In this paper, to overcome the difficulty of the fully-implicit scheme, and reduce the computational complexity for computing the two-dimensional problem, we develop a linear implicit energy-preserving scheme based on the the SAV approach and Fourier pseudo-spectral method to solve the fractional NLS equation.  The proposed scheme is easy to implement and much more efficient for the reason that one only needs to solve  decoupled linear equations with constant coefficients at each time step. Moreover, the Fourier transform technique can be introduced to reduce the computational complexity in the implementation.

This paper is organized as follows. In section 2, we reformulate the equation as a Hamiltonian
system and obtain a new equivalent system via introducing a scalar auxiliary variable. In section 3, first, we use the
Fourier pseudo-spectral method for space discretization of the equivalent system and obtain a semi-discrete conservative scheme. Then utilizing Crank-Nicolson method to discrete the semi-discrete system in time gives a linear implicit scheme in the fully-discrete version. Besides, a fast solver method based on the fast Fourier transformation (FFT) technique is given to practical computation. Some numerical examples are presented in section 4 to illustrate the theoretical results. We draw some conclusions in section 5.

%%%%%%%%%%%%%%%%%%%%%%%%%%%%%%%%%%%%%%%%%% begin section 1 introduction %%%%%%%%%%%%%%%%%%%%%%%%%%%%%%%%%%%%%%%%%%%%
\section{Equivalent system with the SAV approach}

In this section, with the one-dimensional (1D) fractional NLS equation as an example, a new equivalent system is given based on the idea of the SAV approach, and a linear implicit energy-preserving scheme is constructed in next section for the equivalent system.  Similarly, the equivalent system and scheme can be generalized to the two-dimensional (2D) case.

We consider the following 1D fractional NLS equation
\begin{align}\label{NLS:eq:2.1}
\text{i}u_{t}-\gamma(-\Delta)^{\frac{\alpha}{2}}u+(V+\beta|u|^2)u=0,\ x\in \Omega \subset \mathbb{R},\ t \in (0,T],
\end{align}
with the initial condition
\begin{align}\label{NLS:eq:2.2}
u( x,0)=u_0(x),
\end{align}
where $\Omega =(x_{L}, x_{R})$. The fractional Laplacian, under the bounded interval $\Omega =(x_{L}, x_{R})$ with the periodic boundary conditions, can be defined by Fourier series as
\begin{align}\label{NLS:eq:2.3}
-(-\Delta)^{\frac{\alpha}{2}}u(x,t)=-\sum\limits_{k\in \mathbb{Z}}|\upsilon_{k}|^{\alpha}\hat{u}_{k}e^{\text{i}\upsilon_{k}(x-x_{L})},
\end{align}
where $\upsilon_{k}=\frac{2k\pi}{x_{R}-x_{L}}$, $\hat{u}_{k}$ is Fourier coefficient and is given by
\begin{align}\label{NLS:eq:2.4}
\hat{u}_{k}=\frac{1}{x_{R}-x_{L}}\int_{\Omega}u(x,t)e^{-\text{i}\upsilon_{k}(x-x_{L})} dx.
\end{align}

For subsequent theoretical analysis, some lemmas are given as follows.
\begin{lem}\rm{\cite{p34}} For two real periodic functions $\phi$ and $\varphi$, we have
\begin{align}\label{NLS:eq:2.5}
\int_{\Omega}\phi (-\Delta)^{\frac{\alpha}{2}}\varphi dx = \int_{\Omega}\varphi (-\Delta)^{\frac{\alpha}{2}}\phi dx=\int_{\Omega}(-\Delta)^{\frac{\alpha}{4}} \varphi (-\Delta)^{\frac{\alpha}{4}}\phi dx.
\end{align}
\end{lem}

\begin{lem}\rm{\cite{p34}} Assuming the function $G[\psi]$ has the form
\begin{align}\label{NLS:eq:2.6}
G[\psi]=\int_{\Omega}g(\theta, \psi(\theta), (-\Delta)^{\frac{\alpha}{4}}\psi(\theta))d\theta,
\end{align}
where $g$ is a smooth function on the $\Omega$, then the variational derivative of $G(\psi)$ is given as follows
\begin{align}\label{NLS:eq:2.7}
\frac{\delta G}{\delta \psi(\theta)}=\frac{\partial g}{\partial \psi}+(-\Delta)^{\frac{\alpha}{4}}\frac{\partial g}{\partial(-\Delta)^{\frac{\alpha}{4}}\psi}.
\end{align}
\end{lem}

\subsection{Hamiltonian formulation}

By setting $u=p+\text{i}q$, we can rewrite system \eqref{NLS:eq:2.1} as a pair of real-valued equations
\begin{align}\label{NLS:eq:2.8}
p_{t}=(-\Delta)^{\frac{\alpha}{2}}q-Vq-\beta(p^2+q^2)q,
\end{align}
\begin{align}\label{NLS:eq:2.9}
q_{t}=-(-\Delta)^{\frac{\alpha}{2}}p+Vp+\beta(p^2+q^2)p,
\end{align}
with the periodic boundary conditions
\begin{align*}
p(x,t)=p(x+L,t),\ \ \ q(x,t)=q(x+L,t),
\end{align*}
where $L=x_{R}-x_{L}$.

\begin{thm}\rm\label{2SG-lem2.1} The system \eqref{NLS:eq:2.8}-\eqref{NLS:eq:2.9} with the periodic boundary conditions has the following energy and mass conservation laws
\begin{align}\label{NLS:eq:2.10} \mathcal{H}=-\frac{1}{2}\int_{\Omega}\Big[((-\Delta)^{\frac{\alpha}{4}}p)^2+((-\Delta)^{\frac{\alpha}{4}}q)^2-V(p^2+q^2)-\frac{\beta}{2}(p^2+q^2)^2 \Big]dx,
\end{align}
\begin{align}\label{NLS:eq:2.11}
 \mathcal{M}=&\int_{\Omega}(p^2+q^2)dx,
\end{align}
where $\mathcal{H}$ and $ \mathcal{M}$ are energy and mass functional, respectively.
\end{thm}
\begin{prf}\rm By taking the inner products of above system \eqref{NLS:eq:2.8}-\eqref{NLS:eq:2.9} with $q_{t}$ and $-p_t$, respectively, one can deduce that the conservation of energy
\begin{align*}
\frac{d}{d t}\mathcal{H}=0,
\end{align*}
where Lemma 2.1 was used.

Then, computing the inner product of \eqref{NLS:eq:2.8}-\eqref{NLS:eq:2.9} with $p$ and $q$, respectively, one can obtain the conservation of mass
\begin{align*}
\frac{d}{d t}\mathcal{M}=0.
\end{align*}
We finish the proof.
\qed
\end{prf}

Based on the fractional variational  derivative formula in Lemma 2.2, we obtain the following result straightforwardly.

\begin{thm}\rm\label{2SG-lem2.2} The system \eqref{NLS:eq:2.8}-\eqref{NLS:eq:2.9} is an infinite-dimensional canonical Hamiltonian system
\begin{align*}
\left(\begin{array}{c}
		p_t\\
		q_t
		\end{array} \right)=\left(\begin{array}{cc}
		0 & -1 \\
		1 & 0
		\end{array} \right)\left(\begin{array}{c}
		\delta\mathcal{H}/\delta p\\
		\delta\mathcal{H}/\delta q
		\end{array} \right),
\end{align*}
where the energy functional $\mathcal{H}$ is given in Theorem 2.1.
\end{thm}

\subsection{Equivalent system and modified energy conservation law}

In the SAV approach, we introduce a scalar variable $w(t)=\sqrt{\mathcal{E}(t)}$, where
\begin{equation}\label{NLS:eq:2.12}
\mathcal{E}=\frac{1}{4}\int_\Omega\big[\beta(p^2+q^2)^2+2V(p^2+q^2)\big]dx.
\end{equation}
Then, the energy function \eqref{NLS:eq:2.10} is transformed into a modified  formal
\begin{equation}\label{NLS:eq:2.13}
\tilde{\mathcal{H}}=-\frac{1}{2}\int_\Omega \Big[((-\Delta)^{\frac{\alpha}{4}}p)^2+((-\Delta)^{\frac{\alpha}{4}}q)^2\Big]dx+ w^2,
\end{equation} and the corresponding fractional NLS system \eqref{NLS:eq:2.8}-\eqref{NLS:eq:2.9} can be written as
\begin{align}\label{NLS:eq:2.14}
p_t=(-\Delta)^{\frac{\alpha}{2}}q-\mathcal{B}_{2}(p,q) w,
\end{align}
\begin{align}\label{NLS:eq:2.15}
q_t=-(-\Delta)^{\frac{\alpha}{2}}p+ \mathcal{B}_{1}(p,q) w,
\end{align}
\begin{align}\label{NLS:eq:2.16}
w_t=\frac{1}{2}(\mathcal{B}_{1}(p,q), p_{t})+\frac{1}{2}(\mathcal{B}_{2}(p,q), q_{t}),
\end{align}
where $\mathcal{B}_1(p,q)=\frac{\beta(p^2+q^2)p+Vp}{\sqrt{\mathcal{E}(t)}}$, $\mathcal{B}_2(p,q)=\frac{\beta(p^2+q^2)q+Vq}{\sqrt{\mathcal{E}(t)}}$, and $(\cdot, \cdot)$ represents the continuous $L^2$-inner product.

Taking the inner products of \eqref{NLS:eq:2.14}-\eqref{NLS:eq:2.15} with $q_t$, $-p_t$, respectively, and then multiplying \eqref{NLS:eq:2.16} with $2w$, we can deduce a
modified energy conservation law
 \begin{align}\label{NLS:eq:2.17}
 \frac{d}{dt}\Big[-\frac{1}{2}\int_{\Omega}\Big(((-\Delta)^{\frac{\alpha}{4}}p)^2+((-\Delta)^{\frac{\alpha}{4}}q)^2\Big)dx+w^{2}\Big]=0.
 \end{align}
It is observe that equivalent system \eqref{NLS:eq:2.14}-\eqref{NLS:eq:2.16} still preserves energy, but the energy is in terms of the new variables now.

\section{Construction of the energy-preserving scheme}

In this section, we construct a linear implicit energy-preserving scheme for equivalent system \eqref{NLS:eq:2.14}-\eqref{NLS:eq:2.16} of the 1D fractional NLS equation.

\subsection{Structure-preserving spatial discretization}

We choose the mesh size $h:=\frac{x_{R}-x_{L}}{N}$ and the time step $\tau:= \frac{T}{M}$ with integers $N$ and $M$. Denote~$\Omega_{h}=\{x_j |\ x_j=x_{L}+jh, 0\leq j \leq N-1  \}$,
$\Omega_{\tau}=\{t_m|\ t_m=m\tau, 0\leq m \leq M \}$. Let $\mathcal{V}_{h} = \{v|v=(v_0,v_1, \cdots, v_{N-1})^{T}\}$ be the space of grid functions. For a given grid function
$\mathcal{\mathring{V}}_{h} = \{v_{j}^{m}|\ v_{j}^{m}=v(x_j, t_m),  (x_j, t_m)\in \Omega_{h}\times \Omega_{\tau}\}$, we introduce some operators for any mesh function $v_{j}^{m} \in \mathcal{\mathring{V}}_{h} $ as

\begin{align*}
\delta_{t}v_{j}^{m}=\frac{v_{j}^{m+1}-v_{j}^{m}}{\tau},\ v_{j}^{m+1/2}=\frac{v_{j}^{m+1}+v_{j}^{m}}{2},\ \bar{v}_{j}^{m+1/2}=\frac{3{v}_{j}^{m}-{v}_{j}^{m-1}}{2}.
\end{align*}
For any two grid functions $v$,~$s$~$\in$~$\mathcal{V}_{h}$,~we define the discrete inner product as
 \begin{align*}
(v, s)=h\sum\limits_{j=0}^{N-1} v_{j}{\overline{s}}_{j},
 \end{align*} and the discrete maximum norm ($l^{\infty}$-norm) as
\begin{align*}
\|v\|_{{\infty}}=\max\limits_{0 \leq j \leq N-1}|v_j|.
\end{align*}

Let $x_{j}\in \Omega_{h}$ be the Fourier collocation points. Use the interpolation polynomial $I_{N}{v}(x)$ to approximate ${v}_{N}(x)$ of the function ${v}(x)$, and which is defined by
\begin{align*}
I_{N}{v}(x)={v}_{N}(x)=\sum\limits_{k=-N/2}^{N/2}\hat{v}_{k}e^{\text{i}k\mu(x-x_{L})},
\end{align*}
where $\mu=\frac{2\pi}{x_{R}-x_{L}}$, and the coefficient
\begin{align*}
\hat{v}_{k}=\frac{1}{Nc_{k}}\sum\limits_{j=0}^{N-1}{v}(x_{j})e^{-\text{i}k\mu(x-x_{L})},
\end{align*}
where $c_{k}=1$ for $|k|< N/2$, and $c_{k}=2$ for $k=\pm N/2$. Then
 the fractional Laplacian $-(-\Delta)^{\frac{\alpha}{2}}v(x)$ can be approximated by
\begin{align}\label{FKGS:eq:2.8}
-(-\Delta)^{\frac{\alpha}{2}}v_{N}(x_j)=-\sum\limits_{k=-N/2}^{N/2}|k\mu|^{\alpha}\hat{v}(x_{j})e^{-\text{i}k\mu(x_{j}-x_{L})}.
\end{align}
 We denote ${v}_j={v}(x_j)$ and plug $\hat{v}_{k}$ in to \eqref{FKGS:eq:2.8}, then above approximation can be written as a matrix form. To this end,
\begin{align}\label{FKGS:eq:2.10}
-(-\Delta)^{\frac{\alpha}{2}} v_{N}(x_j)&=-\sum\limits_{k=-N/2}^{N/2}|k\mu|^{\alpha}(\frac{1}{Nc_{k}}\sum\limits_{l=0}^{N-1}{v_{l}}e^{-\text{i}k\mu(x-x_{L})})e^{-\text{i}k\mu(x_{j}-x_{L})}\nonumber\\
&=\sum\limits_{l=0}^{N-1}{v}_{l}(-\sum\limits_{k=-N/2}^{N/2}\frac{1}{Nc_{k}}|k\mu|^{\alpha}e^{-\text{i}k\mu(x_j-x_{L})})\nonumber\\
&=(D^{\alpha}{\textbf{v}})_{j},
\end{align}
where ${\textbf{v}}=(v_0, v_1, \cdots, v_{N-1})$, and the spectral differential matrix $D^{\alpha}$ is an $N\times N$ symmetric matrix with the elements \cite{p34}
\begin{align*}
(D^{\alpha})_{j,l}&=-\sum\limits_{k=-N/2}^{N/2}\frac{1}{Nc_{k}}|k\mu|^{\alpha}e^{-\text{i}k\mu(x_j-x_{l})}.
\end{align*}

\begin{rmk} \rm
When $\alpha=2$, the matrix $D^{\alpha}$ will reduce to the second-order spectral differential matrix \cite{p34}.
\end{rmk}

\subsection{Energy-preserving semi-discrete scheme}

%Applying the Fourier pseudo-spectral method in space for the system \eqref{NLS:eq:2.14}-\eqref{NLS:eq:2.16},
We set $P=(P_0, P_1,\cdots, P_{N-1})^{T}$, $Q=(Q_0, Q_1,\cdots, Q_{N-1})^{T}$, $W=(W_0, W_1,\cdots, W_{N-1})^{T}$. Then, discretizing \eqref{NLS:eq:2.14}-\eqref{NLS:eq:2.15} in space by the the Fourier pseudo-spectral method,
 we obtain the following semi-discrete system
\begin{align}\label{NLS:eq:3.1}
\frac{d}{dt}P_{j}=-D^{\alpha}{Q}_{j}-{B}^{j}_{2}W_{j},\end{align}
\begin{align} \label{NLS:eq:3.2}
\frac{d}{dt}Q_{j}=D^{\alpha}{P}_{j}+ {B}^{j}_{1}W_{j},\end{align}
\begin{align}\label{NLS:eq:3.3}
\frac{d}{dt}W_{j}=\frac{1}{2}\Big({B}^{j}_{1}, \frac{d}{dt}P_{j}\Big)+\frac{1}{2}\Big({B}^{j}_{2}, \frac{d}{dt}Q_{j}\Big),
\end{align}
where ${B}^{j}_{1}=\mathcal{B}_{1}(P_j,Q_{j})$, ${B}^{j}_{2}=\mathcal{B}_{2}(P_j,Q_{j})$, $0\leq j\leq N-1$.% ${\textbf{P}}=(P_0, P_1, \cdots, P_{N-1})$, and ${\textbf{Q}}=(Q_0, Q_1, \cdots, Q_{N-1})$.

Next, we will show the energy conservation law of above semi-discrete scheme in the following theorem.

\begin{thm}\rm\label{FSE-lem3.1} The semi-discrete scheme \eqref{NLS:eq:3.1}-\eqref{NLS:eq:3.3} satisfies the energy conservation law
\begin{align*}
{E}(t)={E}(0),
\end{align*}
where
\begin{align*}
{E}(t)=-\frac{1}{2}({ P}^TD^{\alpha}P+ { Q}^TD^{\alpha}Q)+||{W}||^2.
\end{align*}

\end{thm}
\begin{prf}\rm By taking the discrete inner products of \eqref{NLS:eq:3.1}, \eqref{NLS:eq:3.2} with $\frac{d}{dt}Q_{{j}}$, $-\frac{d}{dt}P_{{j}}$, respectively, and then multiplying \eqref{NLS:eq:3.3} with $2W_{{j}}$, one can deduce
\begin{align}\label{NLS:eq:3.4}
(\frac{d}{dt}{P_j},\frac{d}{dt}{Q_j})=(-D^{\alpha}Q_j,\frac{d}{dt}Q_{{j}} )-{h }\sum_{j=0}^{N-1} {B}^{j}_{2}W_{j}\frac{d}{dt} Q_{j},\end{align}
\begin{align}\label{NLS:eq:3.5}
 (\frac{d}{dt}{Q_j},-\frac{d}{dt}{ P_j})=(D^{\alpha}P_j,-\frac{d}{dt}P_{{j}} )-{h}\sum_{j=0}^{N-1} {B}^{j}_{1}W_{j}\frac{d}{dt} P_{j},\end{align}
\begin{align}\label{NLS:eq:3.6}
\frac{d}{dt}||{W_j}||^2=h\sum_{j=0}^{N-1}({B}^{j}_{1}W_{j}\frac{d}{dt}P_{j}+{B}^{j}_{2}W_{j}\frac{d}{dt}Q_{j}).
\end{align}
Substituting \eqref{NLS:eq:3.4} and \eqref{NLS:eq:3.6} into \eqref{NLS:eq:3.5}, we derive
\begin{align}\label{NLS:eq:3.7}
\frac{d}{dt}\Big(-\frac{1}{2}({ P}^TD^{\alpha}P+ { Q}^TD^{\alpha}Q)+\beta||{W}||^2\Big)=0.
\end{align}
Thus, the proof is completed.
\qed
\end{prf}
\subsection{A fully-discrete linear energy-preserving scheme}

Based on the discussions of semi-discrete system \eqref{NLS:eq:3.1}-\eqref{NLS:eq:3.3},
in this subsection, our goal is to establish a linear implicit fully-discrete scheme.

 Applying the Crank-Nicolson method for \eqref{NLS:eq:3.1}-\eqref{NLS:eq:3.3} in time, further utilizing the extrapolation
technique, we can obtain a linear implicit scheme for the fractional NLS equation as follows
\begin{align}\label{NLS:eq:3.8}
\frac{P^{m+1}-P^m}{\tau}=-D^{\alpha}Q^{{m+1/2}}-\bar{B}_2^{m} W^{m+1/2},
\end{align}
\begin{align}\label{NLS:eq:3.9}
\frac{Q^{m+1}-Q^m}{\tau}=D^{\alpha}P^{{m+1/2}}+\bar{B}_1^{m} W^{m+1/2},
\end{align}
 \begin{align}\label{NLS:eq:3.10}
W^{m+1}-W^m=\frac{1}{2}\big(\bar{B}_1^{n},P^{m+1}-P^m\big)+\frac{1}{2}\big(\bar{B}_2^{m},Q^{m+1}-Q^m\big),
\end{align}
where $\bar{B}_1^{m}=\mathcal{B}_1(\bar{P}^{m+1/2},\bar{Q}^{m+1/2})$ and $\bar{B}_2^{m}=\mathcal{B}_1(\bar{P}^{m+1/2},\bar{Q}^{m+1/2})$.

The proposed scheme \eqref{NLS:eq:3.8}-\eqref{NLS:eq:3.10} is known as FSAV scheme. A result on energy of conservation property of the scheme is presented by the
following theorem.

\begin{thm}\rm \label{FSE:thm3.2} The fully-discrete scheme \eqref{NLS:eq:3.8}-\eqref{NLS:eq:3.10} possesses the following discrete total energy conservation law
\begin{align*}
{H}^m={H}^{m+1},\ \ 0\le m\le M-1,
\end{align*}
where
 \begin{align*}
 {H}^m=-\frac{1}{2}\Big((P^{n})^TD^{\alpha}P^{n}+(Q^{m})^TD^{\alpha}Q^{m}\Big)+ \|W^m\|^2.
\end{align*}
\end{thm}
\begin{prf}\rm Taking the inner product of \eqref{NLS:eq:3.8}-\eqref{NLS:eq:3.9} with $(Q^{m+1}-Q^n)/\tau$, $(P^{m+1}-P^n)/\tau$, and multiplying \eqref{NLS:eq:3.10} with $(W^{m+1}+W^m)/\tau$, adding them together, we deduce
\begin{align}\label{NLS:eq:3.11}
\|W^{m+1}\|^2&-\frac{1}{2}\Big((P^{m+1})^TD^{\alpha}P^{m+1}+(Q^{m+1})^TD^{\alpha}Q^{m+1}\Big)\nonumber\\&+\|W^m\|^2-\frac{1}{2}\Big((P^{m})^TDP^{m}-(Q^{m})^TDQ^{m}\Big) =0.
\end{align}
Thus, we have
\begin{align}\label{NLS:eq:3.12}
{H}^m={H}^{m+1}, \ \ 0\le m\le M-1.
\end{align}
This ends the proof.
\qed
\end{prf}

To demonstrate the linear implicit advantage with constant coefficient matrix of the FSAV scheme, we first eliminate $W^{m+1}$ and obtain
\begin{equation}\label{NLS:eq:3.13}
\begin{aligned}
&\frac{P^{m+1}-P^m}{\tau}+D^{\alpha}Q^{{m+1/2}}+\bar{B}_2^{m} W^{m}+\frac{1}{4}\bar{B}_2^{m}\Big[\big(\bar{B}_1^{m},P^{m+1}-P^m\big)+\big(\bar{B}_2^{m},Q^{m+1}-Q^m\big)\Big]=0,
\end{aligned}
\end{equation}
\begin{equation}\label{NLS:eq:3.14}
\begin{aligned}
&\frac{Q^{m+1}-Q^m}{\tau}-D^{\alpha}P^{{m+1/2}}-\bar{B}_1^{m} W^{m}-\frac{1}{4}\bar{B}_1^{m}\Big[\big(\bar{B}_1^{m},P^{m+1}-P^m\big)+\big(\bar{B}_2^{m},Q^{n+1}-Q^m\big)\Big]=0,
\end{aligned}
\end{equation}
which can be further arranged as
\begin{equation}\label{NLS:eq:3.15}
\begin{aligned}
&P^{m+1}+\frac{\tau}{2} D^{\alpha}Q^{{m+1}}+\frac{\tau}{4}\bar{B}_2^{m}\Big[\big(\bar{B}_1^{m},P^{m+1}\big)+\big(\bar{B}_2^{m},Q^{m+1}\big)\Big]\\
&\qquad\qquad=P^n-\frac{\tau}{2} D^{\alpha}Q^{{m}}-\tau\bar{B}_2^{m} W^{m}+\frac{\tau}{4}\bar{B}_2^{m}\Big[\big(\bar{B}_1^{m},P^m\big)+\big(\bar{B}_2^{m},Q^m\big)\Big]:=C_1^m,
\end{aligned}
\end{equation}
\begin{equation}\label{NLS:eq:3.16}
\begin{aligned}
&Q^{m+1}-\frac{\tau}{2} D^{\alpha}P^{{m+1}}-\frac{\tau }{4}\bar{B}_1^{m}\Big[\big(\bar{B}_1^{m},P^{m+1}\big)+\big(\bar{B}_2^{m},Q^{m+1}\big)\Big]\\
&\qquad\qquad=Q^m+\frac{\tau}{2} D^{\alpha}P^{{m}}+\tau \bar{B}_1^{m} W^{m}-\frac{\tau }{4}\bar{B}_1^{m}\Big[\big(\bar{B}_1^{m},P^m\big)+\big(\bar{B}_2^{m},Q^m\big)\Big]:=C_2^m.
\end{aligned}
\end{equation}
Let $Z^m=\left(\begin{array}{c}
P^m\\
Q^m
\end{array} \right)$, $\bar{B}^m=\left(\begin{array}{c}
\bar{B}_1^m\\
\bar{B}_2^m
\end{array} \right)$, $\bar{B}_{r}^m=\left(\begin{array}{c}
\bar{B}_2^m\\
-\bar{B}_1^m
\end{array} \right)$, ${C}^m=\left(\begin{array}{c}
{C}_1^m\\
{C}_2^m
\end{array} \right)$, and denote $A=\left(\begin{array}{cc}
	I & \frac{\tau}{2}D^{\alpha} \\
	-\frac{\tau}{2}D^{\alpha} &  I
\end{array} \right)$, we can derive a compact scheme for the above system, which reads
\begin{equation}\label{NLS:eq:3.17}
AZ^{m+1}+\frac{\tau}{4}\big(\bar{B}^{m},Z^{m+1}\big)\bar{B}_{r}^m=C^m.
\end{equation}
Multiplying \eqref{NLS:eq:3.17} with $A^{-1}$, and taking the discrete
inner product with $\bar{B}^m$, we can obtain
\begin{equation}\label{NLS:eq:3.18}
\big(\bar{B}^{m},Z^{m+1}\big)=\dfrac{\big(\bar{B}^m,A^{-1}C^m\big)}{1+\frac{\tau}{4}\chi^m},
\end{equation}
where $\chi^m=(\bar{B}^m,A^{-1}\bar{B}_r^m)$. Finally, substituting the result into \eqref{NLS:eq:3.17} will give the solution of $Z^{m+1}$.  To summarize, the FSAV scheme \eqref{NLS:eq:3.8}-\eqref{NLS:eq:3.10} can be easily implemented in the following manner:

(i) First, compute $\chi^m=(\bar{B}^m,A^{-1}\bar{B}_r^m)$. This can be completed by solving a algebraic system, which has constant coefficients.

(ii) Second, compute $\big(\bar{B}^{m},Z^{m+1}\big)$ using \eqref{NLS:eq:3.18}. This requires solving another system $ A^{-1}C^{m}$ with constant coefficients.

(iii) Finally, with $\big(\bar{B}^{m},Z^{m+1}\big)$, $(\bar{B}^m,A^{-1}\bar{B}_r^m)$ and $ A^{-1}C^{m}$ known, we can obtain $Z^{m+1}$ from \eqref{NLS:eq:3.17}.

In summary, we only need to solve two linear systems of \eqref{NLS:eq:3.17} and \eqref{NLS:eq:3.18} successively with constant matrix $A$.

Next, we give a fast solver for the computation of $A^{-1}C^m$ or $A^{-1}B_r^m$. The spectral differential matrix $D^{\alpha}$ that can be further written as $D^{\alpha}$=$F^{-1}\Lambda F$, where $F$ and $F^{-1}$ are the corresponding matrices for discrete Fourier transformation, $\Lambda$ is a diagonal matrix with eigenvalues of $D^{\alpha}$ being its entries. Then we can rewrite the coefficient matrix $A$ as
\begin{equation}\label{NLS:eq:3.19}
A=\mathbb{F}^{-1}\mathbb{M}\mathbb{F},
\end{equation}
where
\begin{equation}\label{NLS:eq:3.20}
\mathbb{F}=\left(\begin{array}{cc}
F & 0 \\
0 &  F
\end{array} \right),\quad \mathbb{M}=\left(\begin{array}{cc}
I & \frac{\tau}{2}\Lambda \\
-\frac{\tau}{2}\Lambda &  I
\end{array} \right).
\end{equation}
Moreover, the inverse of $A$ satisfies
\begin{equation}\label{NLS:eq:3.21}
A^{-1}=\mathbb{F}^{-1}\mathbb{M}^{-1}\mathbb{F},
\end{equation}
where $\mathbb{M}^{-1}$ can be calculated explicitly as
\[
\mathbb{M}^{-1}=\mathbb{M}^T\left(\begin{array}{cc}
K & 0 \\
0 &  K
\end{array} \right),\quad \mbox{with}\quad K=(I+\frac{\tau^2}{4}\Lambda^2)^{-1}.
\]
Notice that $\Lambda$ is a diagonal matrix and therefore the computational cost of its inverse is negative.

In the 2D case, \eqref{NLS:eq:3.15}-\eqref{NLS:eq:3.16} becomes
\begin{equation}\label{NLS:eq:3.22}
\begin{aligned}
P^{m+1}&+\frac{\tau}{2} \big(D^{\alpha}_xQ^{{m+1}}+Q^{{m+1}}(D^{\alpha}_y)^T\big)+\frac{\tau}{4}\bar{B}_2^{m}\Big[\big(\bar{B}_1^{m},P^{m+1}\big)+\big(\bar{B}_2^{m},Q^{m+1}\big)\Big]
\\&=P^m-\frac{\tau}{2} D_x^{\alpha}Q^{{m}}-\frac{\tau}{2} Q^{{m}}(D_{y}^{\alpha})^{T}-\tau\bar{B}_2^{m} W^{m}+\frac{\tau}{4}\bar{B}_2^{m}\Big[\big(\bar{B}_1^{m},P^m\big)+\big(\bar{B}_2^{m},Q^m\big)\Big]:=C_1^m,
\end{aligned}
\end{equation}
\begin{equation}\label{NLS:eq:3.23}
\begin{aligned}
Q^{m+1}&-\frac{\tau}{2} \big(D^{\alpha}_xQ^{{m+1}}+Q^{{m+1}}(D^{\alpha}_y)^T\big)-\frac{\tau }{4}\bar{B}_1^{m}\Big[\big(\bar{B}_1^{m},P^{m+1}\big)+\big(\bar{B}_2^{m},Q^{m+1}\big)\Big]
\\&=Q^m+\frac{\tau}{2} D_x^{\alpha}P^{{m}}+\frac{\tau}{2} P^{{m}}(D_{y}^{\alpha})^{T}+\tau\bar{B}_1^{m} W^{m}-\frac{\tau}{4}\bar{B}_1^{m}\Big[\big(\bar{B}_1^{m},P^m\big)+\big(\bar{B}_2^{m},Q^m\big)\Big]:=C_2^m,
\end{aligned}
\end{equation}
where $D^{\alpha}_x$ and $D^{\alpha}_y$ are the corresponding spectral differential matrices with respect to $x$ and $y$ directions and the variables are now of matrix forms. Firstly, we transform the scheme of unknowns $P^{m+1}$, $Q^{m+1}$ from matrix forms to vector forms $\bm p^{m+1}$, $\bm q^{m+1}$, and follow the above strategies \eqref{NLS:eq:3.15}-\eqref{NLS:eq:3.18}, which leads to solve twice the algebraic system $\bm A\bm x=\bm b$ with
\[
\begin{aligned}
\bm A&=\left(\begin{array}{cc}
I_y\otimes I_x &  \frac{\tau}{2}(I_y\otimes D^{\alpha}_x+D^{\alpha}_y\otimes I_x)\\ [1ex]
-\frac{\tau}{2}(I_y\otimes D^{\alpha}_x+D^{\alpha}_y\otimes I_x) &  I_y\otimes I_x
\end{array} \right)\\
%&=\left(\begin{array}{cc}
%F_y^{-1}F_y\otimes F_x^{-1}F_x &  F_y^{-1}F_y\otimes F_x^{-1}\Lambda_xF_x+F_x^{-1}\Lambda_yF_x\otimes F_x^{-1}F_x\\  [1ex]
%-F_y^{-1}F_y\otimes F_x^{-1}\Lambda_xF_x+F_y^{-1}\Lambda_yF_y\otimes F_x^{-1}F_x &  F_y^{-1}F_y\otimes F_x^{-1}F_x
%\end{array} \right)\\[1ex]
%&=\left(\begin{array}{cc}
%(F_y^{-1}\otimes F_x^{-1})(F_y\otimes F_x) &  (F_y^{-1}\otimes F_x^{-1})(I_y\otimes \Lambda_x+\Lambda_y\otimes I_x)(F_y\otimes F_x)\\ [1ex]
%-(F_y^{-1}\otimes F_x^{-1})(I_y\otimes \Lambda_x+\Lambda_y\otimes I_x)(F_y\otimes F_x) &  (F_y^{-1}\otimes F_x^{-1})(F_y\otimes F_x)
%\end{array} \right)\\[1ex]
&=\mathbb{F}^{-1}\mathbb{M}\mathbb{F},
\end{aligned}
\]
where
\[
\mathbb{F}=\left(\begin{array}{cc}
F_y\otimes F_x &  0\\ [1ex]
0 &  F_y\otimes F_x
\end{array} \right),\quad \mathbb{M}=\left(\begin{array}{cc}
I_y\otimes I_x &  \frac{\tau}{2}(I_y\otimes \Lambda_x+\Lambda_y\otimes I_x)\\ [1ex]
-\frac{\tau}{2}(I_y\otimes \Lambda_x+\Lambda_y\otimes I_x) &  I_y\otimes I_x
\end{array} \right),
\]
and $I_x$, $I_y$ are identical matrices, $\Lambda_x$, $\Lambda_y$ correspond to eigenvalues of spectral differential matrices $D^{\alpha}_x$, $D^{\alpha}_y$, respectively.
The inverse of $\bm A$ can then be calculated by
$\mathbb{F}^{-1}\mathbb{M}^{-1}\mathbb{F}$,
where
\[
\mathbb{M}^{-1}=\mathbb{M}^T\left(\begin{array}{cc}
K & 0 \\
0 &  K
\end{array} \right),\quad \mbox{with}\quad K=\big(I_y\otimes I_x+\frac{\tau^2}{4}(I_y\otimes \Lambda_x+\Lambda_y\otimes I_x)^2\big)^{-1}.
\]
Also $K$ is just the inverse of a diagonal matrix and can be efficiently obtained.

\begin{rmk} \rm  Based on the Remark 2.1, when $\alpha=2$, the FSAV scheme is suitable for the classical Schr\"{o}dinger equation.
\end{rmk}

\section{Numerical examples}

In this section,  numerical examples are presented to show the accuracy and energy conservation property of FSAV scheme \eqref{NLS:eq:3.8}-\eqref{NLS:eq:3.10} and
simulate the fractional NLS equation with  potential function.

To obtain numerical errors, we use the error function defined as follows
\begin{align*}
E(\tau)=\|P^{M}_{N}-P^{2M}_{N}\|_{\infty}+\|Q^{M}_{N}-Q^{2M}_{N}\|_{\infty},
\end{align*}
\begin{align*}
E(N)=\|P^{M}_{N}-P^{M}_{2N}\|_{\infty}+\|Q^{M}_{N}-Q^{M}_{2N}\|_{\infty},
\end{align*}
where $\|P_{N}^{M}-P_{N}^{2M}\|_{\infty}:=\|P(\frac{T}{M},\frac{L}{N})-P(\frac{T}{2M},\frac{L}{N})\|_{\infty}$, $\|P^{M}_{N}-P^{M}_{2N}\|_{\infty}:=\|P(\frac{T}{M},\frac{L}{N})-P(\frac{T}{M},\frac{L}{2N})\|_{\infty}$, etc.,
and the convergence orders in time and space of the $l^{\infty}$-norm errors on two successive time step sizes $\tau$ and $\tau/2$ and two successive polynomial degrees $N$ and $2N$ are calculated as
\begin{align*}
\text{order}=\left\{\begin{array}{lll}
              &{\text{log}_{2}{[{E}(\tau)/{E}(\tau/2)]}},\ &\text{ in time}, \\
 \\
&{\text{log}_{2}[{E}( N)/{E}( 2N)]},\ &\text{ in space}.
\\
             \end{array}
\right.
\end{align*}
The relative errors of energy and mass
are defined as
\begin{align*}
RH^{m}=|(H^{m}-H^{0})/H^{0}|,\ \ \ RM^{m}=|(M^{m}-M^{0})/M^{0}|,
\end{align*}
where $H^{m}$ and $M^{m}$ denote the energy and mass at $t=m\tau$, respectively.

\subsection{Simulations of the 1D fractional NLS equation}

 In this subsection, we study the 1D fractional NLS equation \eqref{NLS:eq:2.1} without potential function as follows
 \begin{align}\label{NLS:eq:4.1}
\text{i}u_{t}-(-\Delta)^{\frac{\alpha}{2}}u+2|u|^2u=0,
\end{align}
with the initial condition
\begin{align}\label{NLS:eq:4.2}
u(x,0)=\text{exp}(-x^2)\exp(-\text{i}x).
\end{align}

First of all, we set the space interval $\Omega=[-16,16]$ to check the accuracy and efficiency of the FSAV scheme. Without loss of generality, we take $\alpha=1.4, 1.7, 1.9, 2$.
Table. 1 shows the temporal errors and convergence orders with $N=256$ at $T=1$. As illustrated in Table. 1, the proposed scheme for the fractional NLS equation for different $\alpha$ is second-order of convergence in time. We show the space errors and convergence orders in Table. 2 when $\tau=10^{-6}$, which shows that for the
sufficiently smooth problem, our scheme is of spectral accuracy in spatial.

\begin{table}[H]
	\centering
	\caption{The numerical errors and convergence orders in time for different $\alpha$ with $N=256$ at $T=1$.}\label{tab:5}
	\begin{tabular*}{1\textwidth}[h]{@{\extracolsep{\fill}}lllllllll} \hline
		\multirow{2}{*}{$\tau$}& \multicolumn{2}{c}{$\alpha=1.4$} &\multicolumn{2}{c}{$\alpha=1.7$} & \multicolumn{2}{c}{$\alpha=1.9$} &  \multicolumn{2}{c}{$\alpha=2.0$}\\ \cline{2-3}\cline{4-5}\cline{6-7}\cline{8-9}
		& error & order & error & order & error & order & error & order\\ \hline
		$0.01$    &  5.65e-05  &  -   &  5.34e-05 &    -  &  6.62e-05 & -    & 7.79e-05 & -\\[1ex]
		$0.005$   &  1.41e-05  & 2.00 &  1.34e-05 & 1.99  &  1.66e-05 & 1.99 & 1.92e-05 &  2.00\\[1ex]
		$0.0025$  &  3.53e-06  & 2.00 &  3.36e-06 & 2.00  &  4.16e-06 & 2.00 & 4.82e-06 &  2.00 \\ [1ex] 	
		$0.00125$ &  8.83e-07  & 2.00 &  8.41e-07 & 2.00  &  1.04e-06 & 2.00 & 1.20e-06 &  2.00\\ [1ex] \hline			
	\end{tabular*}
\end{table}

\begin{table}[H]
	\centering
	\caption{The spatial errors and convergent orders for different $\alpha$ with $\tau= 10^{-6}$ at $T=1$.}\label{tab:5}
	\begin{tabular*}{1\textwidth}[h]{@{\extracolsep{\fill}}lllllllll} \hline
		\multirow{2}{*}{$N$}& \multicolumn{2}{c}{$\alpha=1.4$} &\multicolumn{2}{c}{$\alpha=1.7$} & \multicolumn{2}{c}{$\alpha=1.9$} &  \multicolumn{2}{c}{$\alpha=2.0$}\\ \cline{2-3}\cline{4-5}\cline{6-7}\cline{8-9}
		& error & order & error & order & error & order & error & order\\ \hline
		$32$    & 3.77e-01  &  -    &  2.02e-01 & -     &  1.46e-01  & -      & 1.28e-01   & -   \\[1ex]
		$64$    & 5.33e-02  & 2.82  &  1.23e-02 & 4.03  &  5.84e-03  & 4.64   & 3.62e-03   &  5.15\\[1ex]
		$128$   & 6.89e-04  & 6.28  &  1.09e-05 & 10.13 &  2.15e-06  & 11.40  & 1.12e-06   &  11.65 \\ [1ex] 	
		$256$   & 1.35e-08  & 15.63 &  4.34e-11 & 17.94 &  3.51e-11  & 15.90  & 3.65e-11   &  14.90\\ [1ex] \hline			
	\end{tabular*}
\end{table}

In section 3, we demonstrate the FSAV scheme is not only linearly implicit, but also induce algebraic systems with constant coefficient matrix and which can be efficiently solved. While the energy-preserving Crank-Nicolson Fourier pseudo-spectral (CNF) method is fully-implicit and needs nonlinear iterations to obtain the numerical solution, which can approximately represent the general implicit method, regarding the computational cost. In Fig. 1, we present comparisons on the computational costs of the two schemes with $T=100$. Obviously, we can conclude that the FSAV scheme significantly reduces the computational cost. Therefore, it is preferable to construct linearly
implicit schemes through the SAV approach for large scale simulations, keeping the system energy being
preserved as well.

\begin{figure}[H]
\centering\begin{minipage}[t]{70mm}
\includegraphics[width=70mm]{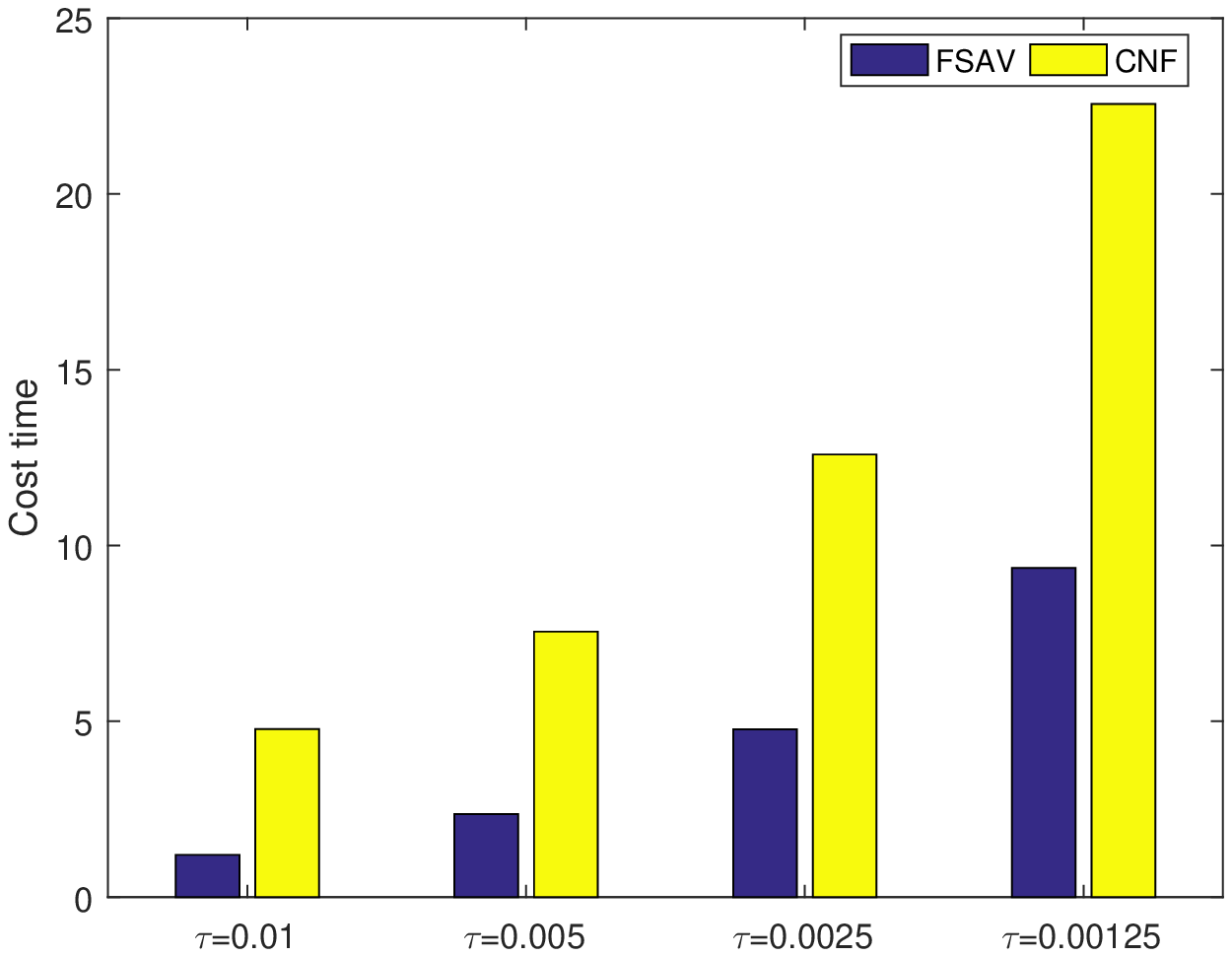}\\
{\footnotesize  \centerline {(a) $\alpha=1.7$}}
\end{minipage}
\begin{minipage}[t]{70mm}
\includegraphics[width=70mm]{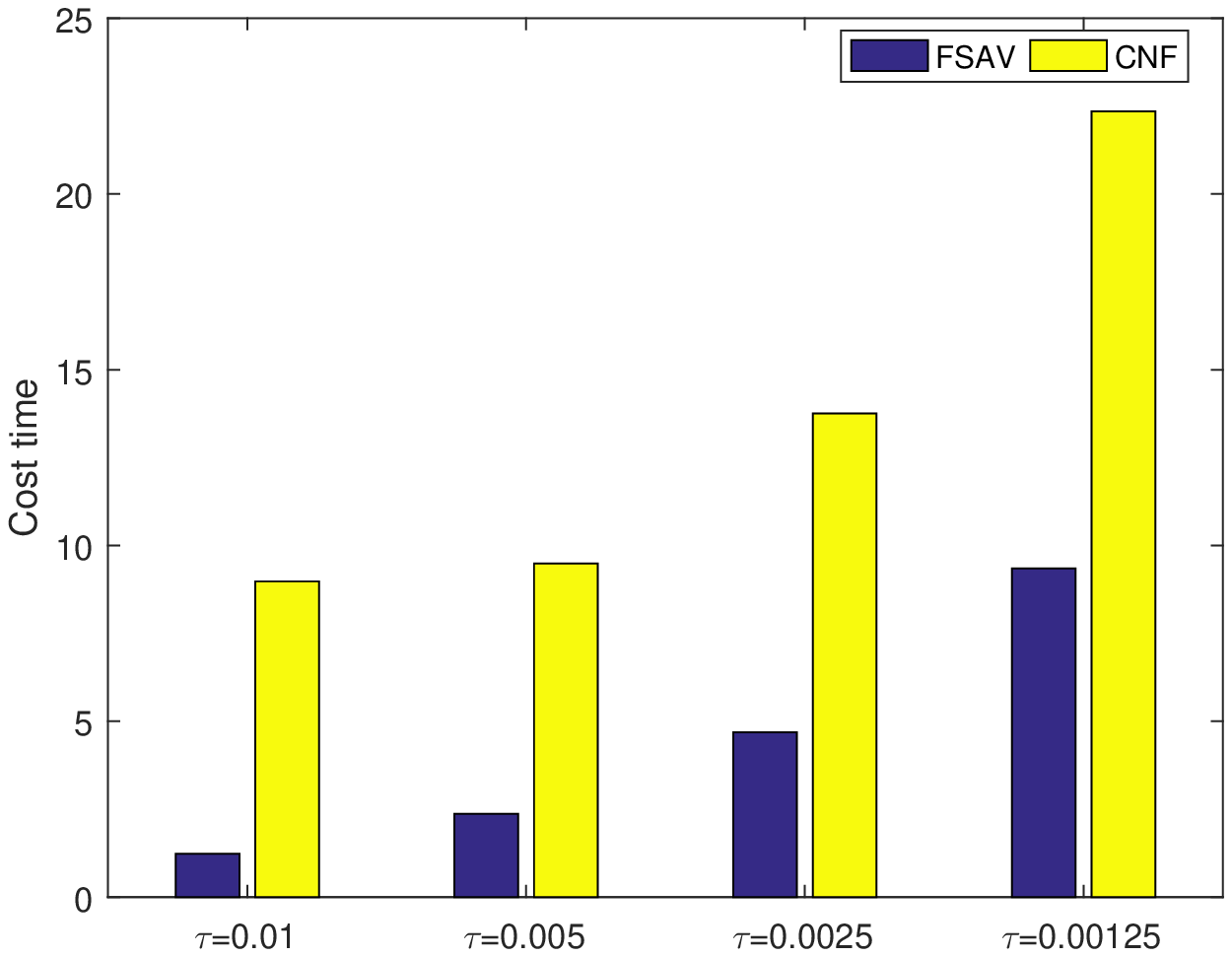}\\
{\footnotesize  \centerline {(b) $\alpha=2$}}
\end{minipage}
\caption{\small {\small{Computational cost of the FSAV scheme and the CNF scheme with different time step.} }}\label{fig522}
\end{figure}

Follow by, we enlarge the computational domain $\Omega=[-40,40]$ to test the discrete conservation law of the FSAV scheme.
In Fig. 2, we depict the the relative errors of the energy $H$ and mass $M$ for different values of fractional order $\alpha$.
It is observed that the FSAV scheme preserves the energy very well in discrete sense, but can not preserve the discrete mass.

\begin{figure}[H]
\centering\begin{minipage}[t]{70mm}
\includegraphics[width=70mm]{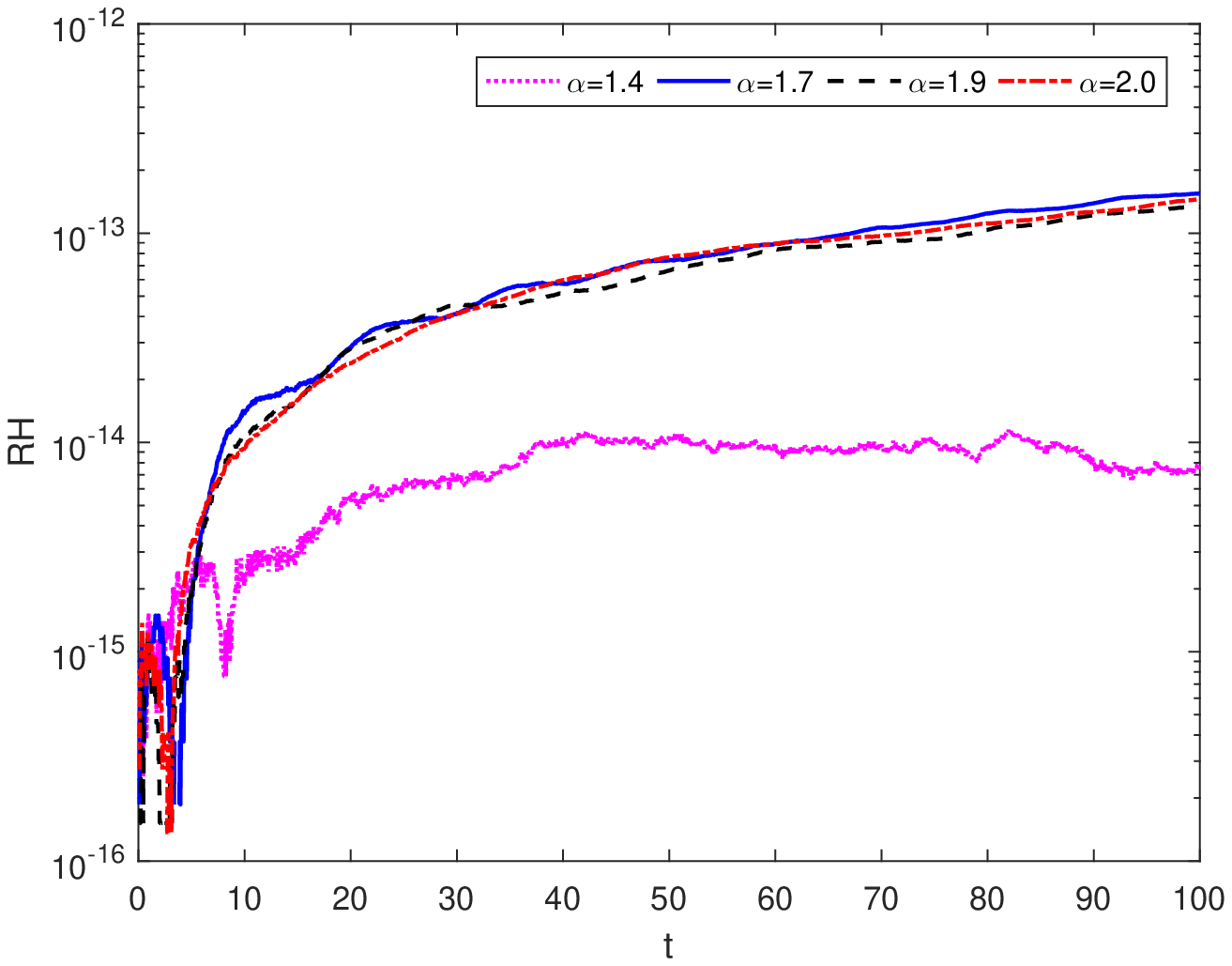}\\
{\footnotesize  \centerline {(a) Relative errors of energy}}
\end{minipage}
\begin{minipage}[t]{70mm}
\includegraphics[width=70mm]{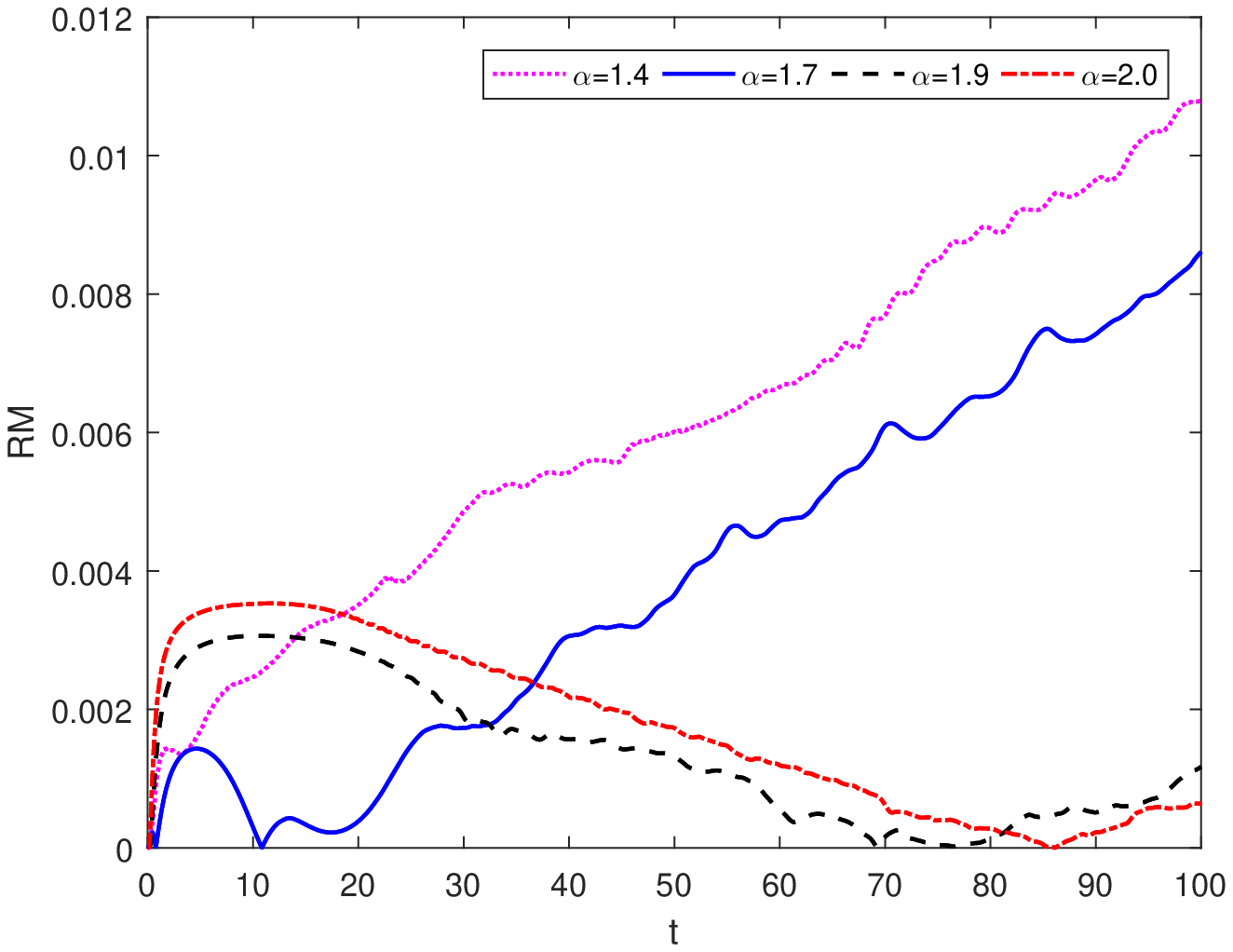}\\
{\footnotesize  \centerline {(b) Relative errors of mass}}
\end{minipage}
\caption{\small {\small{Relative errors of energy and mass with $h=0.5, \tau=0.01$ for different $\alpha$.} }}\label{fig522}
\end{figure}

Finally, we pour attention into the relationship between the fractional order $\alpha$ and the shape of the solition.
Fig. 3 presents the numerical results for different values of $\alpha$. We can find that the fractional order $\alpha$ will affect the shape of the soliton, the shape of the soliton changes dramatically when
the fractional order $\alpha$ changes from 1.4 to 2. The results indicate that the solitary wave can propagate in a stable way under finite initial conditions even after the propagation over sufficiently long times.

\begin{figure}[H]
\centering\begin{minipage}[t]{70mm}
\includegraphics[width=70mm]{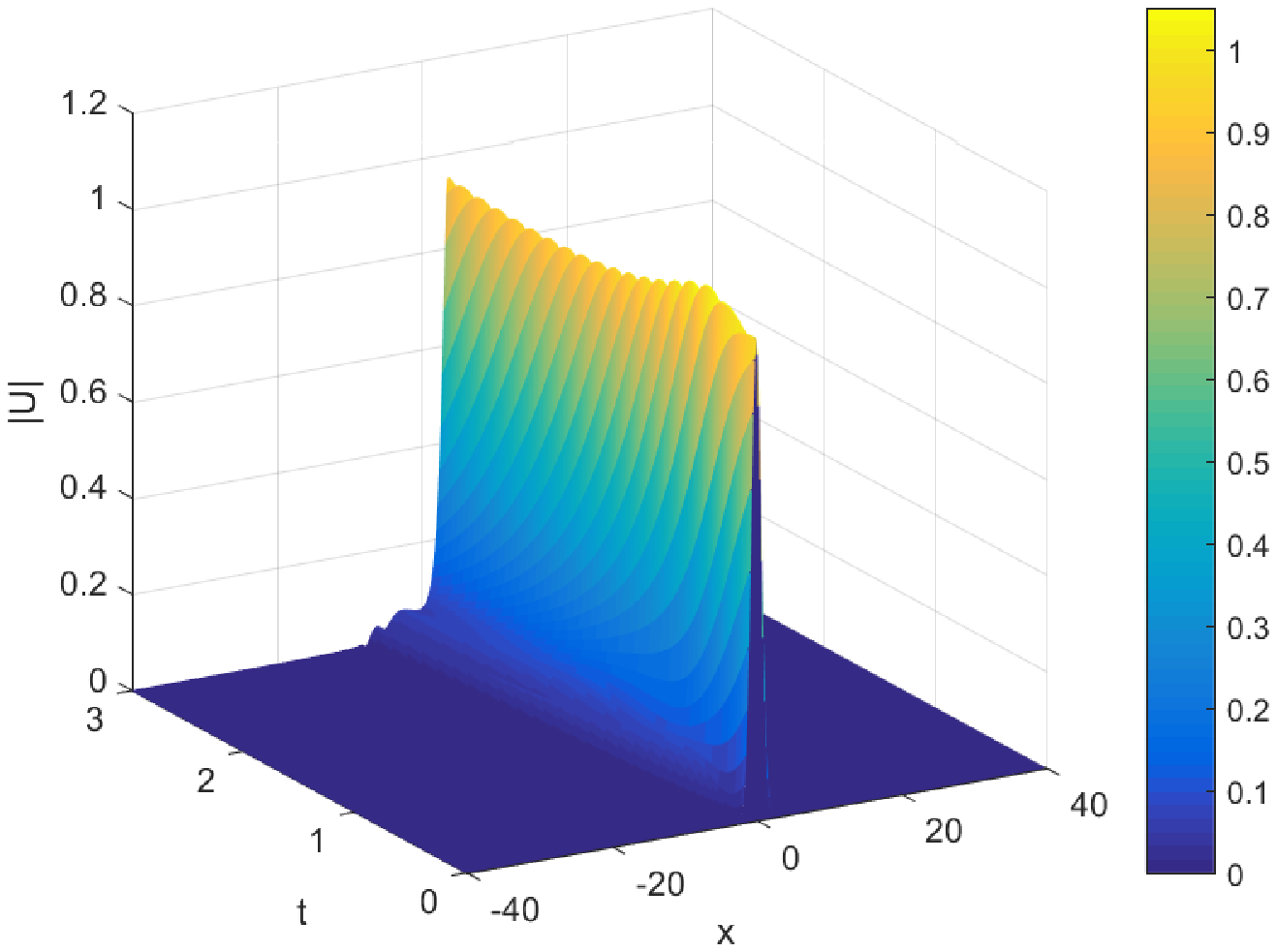}\\
{\footnotesize  \centerline {(a) $\alpha=1.4$}}
\end{minipage}
\begin{minipage}[t]{70mm}
\includegraphics[width=70mm]{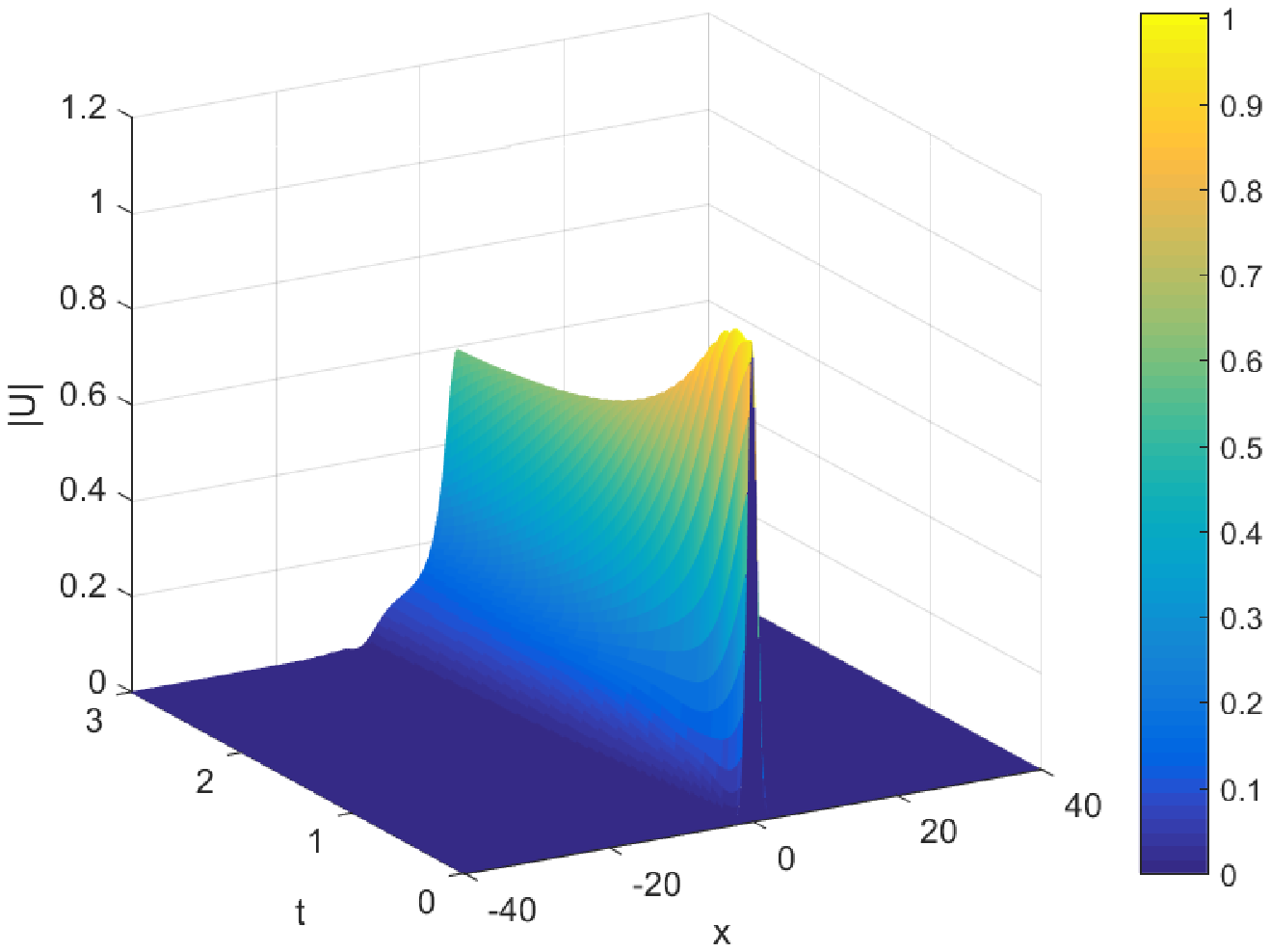}\\
{\footnotesize  \centerline {(b) $\alpha=1.7$}}
\end{minipage}
\begin{minipage}[t]{70mm}
\includegraphics[width=70mm]{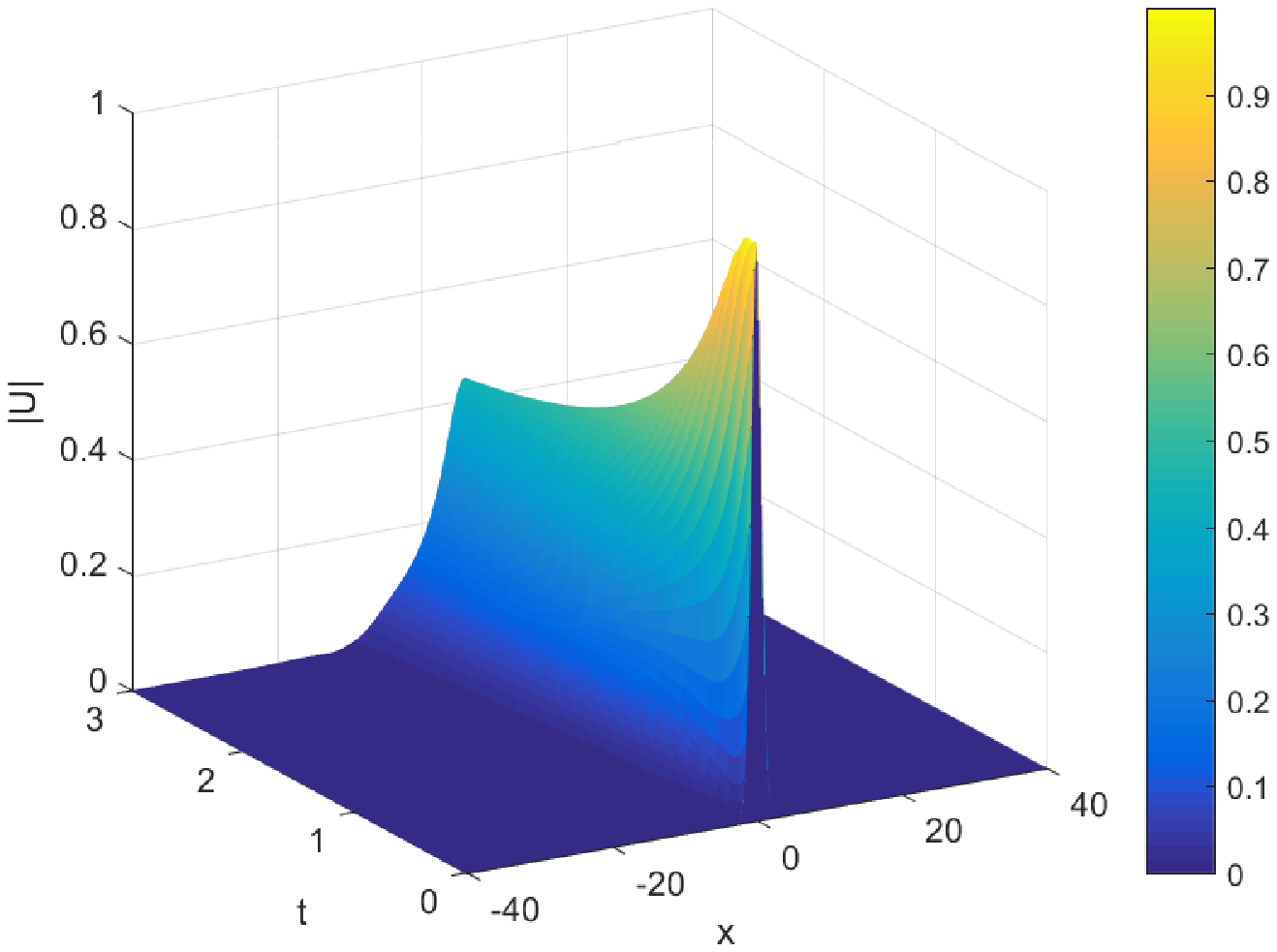}\\
{\footnotesize  \centerline {(c) $\alpha=1.9$}}
\end{minipage}
\begin{minipage}[t]{70mm}
\includegraphics[width=70mm]{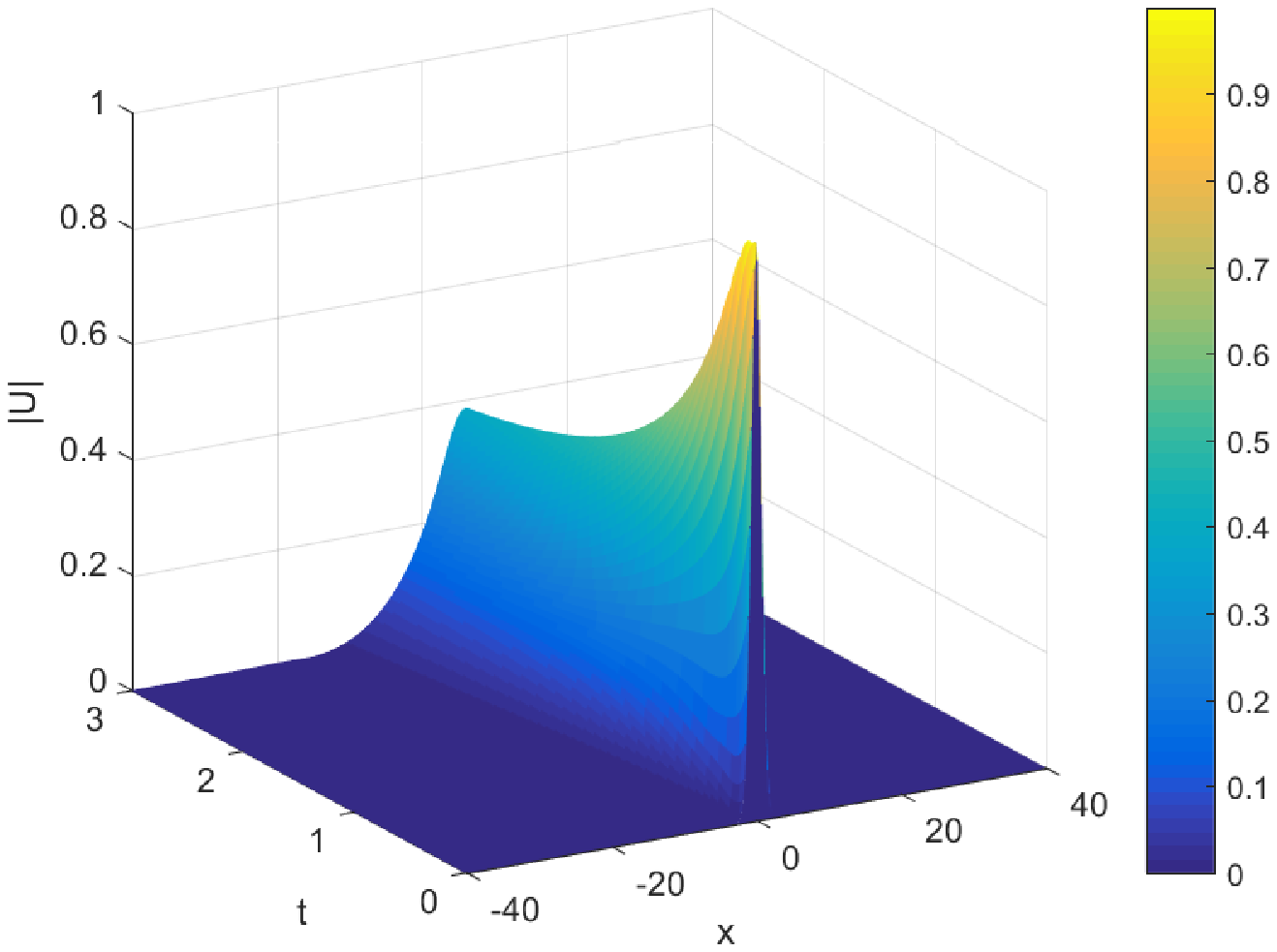}\\
{\footnotesize   \centerline {(d) $\alpha=2.0$}}
\end{minipage}
\caption{\small {\small{Evolution of the solitons with $N=400, \tau=0.001$ for different order $\alpha$.} }}\label{fig522}
\end{figure}

\subsection{Simulations of the 2D fractional NLS equation}

In this subsection, we use the same number of points $N$ in the $x$ and $y$ directions
and set $h=h_x=h_y=\frac{x_{R}-x_{L}}{N}$.

\textbf{Example 4.2.} We consider the 2D fractional NLS equation (\ref{NLS:eq:1.1}) with $\gamma=\beta=1$ and $V=0$. The initial condition is chosen as
\begin{align}\label{NLS:eq:4.3}
u(x,y,0)=\frac{2}{\sqrt{\pi}}\text{exp}(-x^2-y^2).
\end{align}
Here, we take different fractional order, i.e., $\alpha=1.3, 1.6,1.9,2$.

First, we test the accuracy of the FSAV scheme. In our computation, we set the space interval $\Omega=[-8,8]\times[-8,8]$.
Table. 3 shows the temporal errors and convergence orders with $N=128$ at $T=1$. The convergence orders
indicate that our scheme is of second-order accuracy in time. Next, we fix the time step $\tau= 10^{-4}$
to compute the space errors of the proposed scheme. The spatial errors for the 2D fractional NLS equation with different $\alpha$ are presented in Table. 4. One can easy see that the spatial errors are very small. It confirms that for the sufficiently smooth problem, the FSAV scheme is of spectral accuracy in spatial direction.

\begin{table}[H]
	\centering
	\caption{The time errors and convergence orders for different $\alpha$ with $N=128$ at $T=1$.}\label{tab:5}
	\begin{tabular*}{1\textwidth}[h]{@{\extracolsep{\fill}}lllllllll} \hline
		\multirow{2}{*}{$\tau$}& \multicolumn{2}{c}{$\alpha=1.3$} &\multicolumn{2}{c}{$\alpha=1.6$} & \multicolumn{2}{c}{$\alpha=1.9$} &  \multicolumn{2}{c}{$\alpha=2.0$}\\ \cline{2-3}\cline{4-5}\cline{6-7}\cline{8-9}
		& error & order & error & order & error & order & error & order\\ \hline
		$0.02$   &  5.80e-04  &  -   &  5.44e-04 &    -  &  6.83e-04 & -    & 1.10e-03 & -\\[1ex]
		$0.01$   &  1.46e-04  & 1.99 &  1.43e-04 & 1.92  &  1.76e-04 & 1.95 & 2.65e-04 &  2.05\\[1ex]
		$0.005$  &  3.67e-05  & 1.99 &  3.62e-05 & 1.98  &  4.41e-05 & 2.00 & 6.02e-05 &  2.13 \\ [1ex] 	
		$0.0025$ &  9.21e-06  & 2.00 &  9.10e-06 & 1.99  &  1.10e-05 & 2.00 & 1.47e-05 &  2.03\\ [1ex] \hline			
	\end{tabular*}
\end{table}

\begin{table}[H]
	\centering
	\caption{The spatial errors and convergent orders for different $\alpha$ with $\tau=10^{-4}$ at $T=1$.}\label{tab:5}
	\begin{tabular*}{1\textwidth}[h]{@{\extracolsep{\fill}}lllllllll} \hline
		\multirow{2}{*}{$N$}& \multicolumn{2}{c}{$\alpha=1.3$} &\multicolumn{2}{c}{$\alpha=1.6$} & \multicolumn{2}{c}{$\alpha=1.9$} &  \multicolumn{2}{c}{$\alpha=2.0$}\\ \cline{2-3}\cline{4-5}\cline{6-7}\cline{8-9}
		& error & order & error & order & error & order & error & order\\ \hline
		$16$  & 8.25e-01  &  -    &  5.27e-01 & -     &  4.75e-01  & -      & 4.73e-01   & -   \\[1ex]
		$32$  & 9.70e-03  & 6.41  &  2.68e-03 & 7.61  &  1.04e-03  & 8.83   & 7.98e-04   &  9.21\\[1ex]
		$64$  & 4.57e-05  & 7.72  &  4.71e-06 & 9.15  &  6.88e-07  & 10.56  & 4.11e-07   &  10.92 \\ [1ex] 	
		$128$ & 1.71e-09  & 14.70 &  1.07e-11 & 18.74 &  5.30e-13  & 20.30  & 5.07e-13   &  19.63\\ [1ex] \hline			
	\end{tabular*}
\end{table}

Second, we enlarge the computational domain $\Omega=[-10,10]\times[-10,10]$ and study the conservative property of the proposed scheme. Here, we take $h=0.5,$~$\tau=0.02$ and compute the discrete energy and mass. Fig. 4 shows the relative errors of energy and mass for different values of fractional order $\alpha$. Numerical results demonstrate that the FSAV scheme only can preserve the energy exactly in discrete sense.

\begin{figure}[H]
\centering\begin{minipage}[t]{70mm}
\includegraphics[width=70mm]{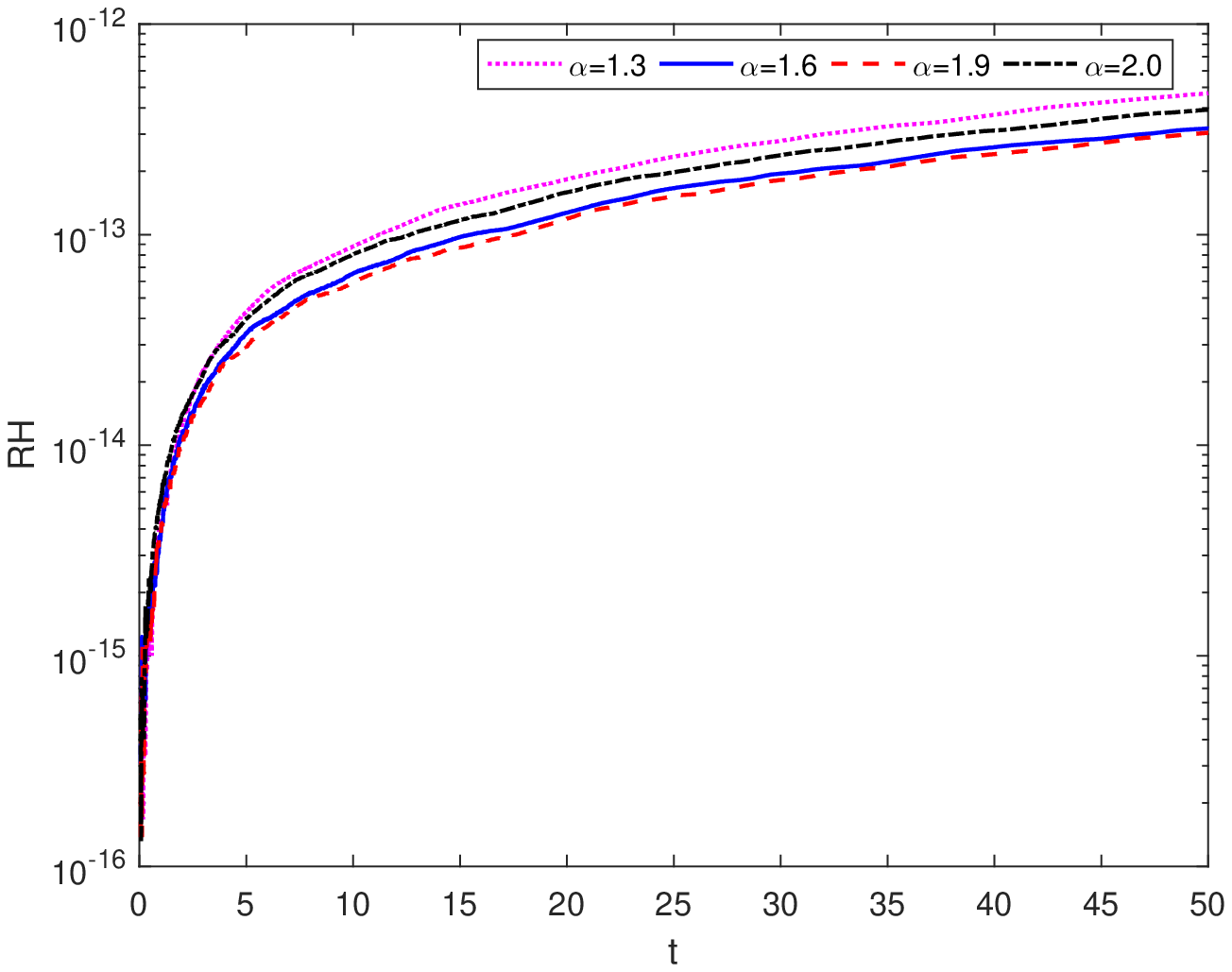}\\
{\footnotesize  \centerline {(a) Relative errors of energy}}
\end{minipage}
\begin{minipage}[t]{70mm}
\includegraphics[width=70mm]{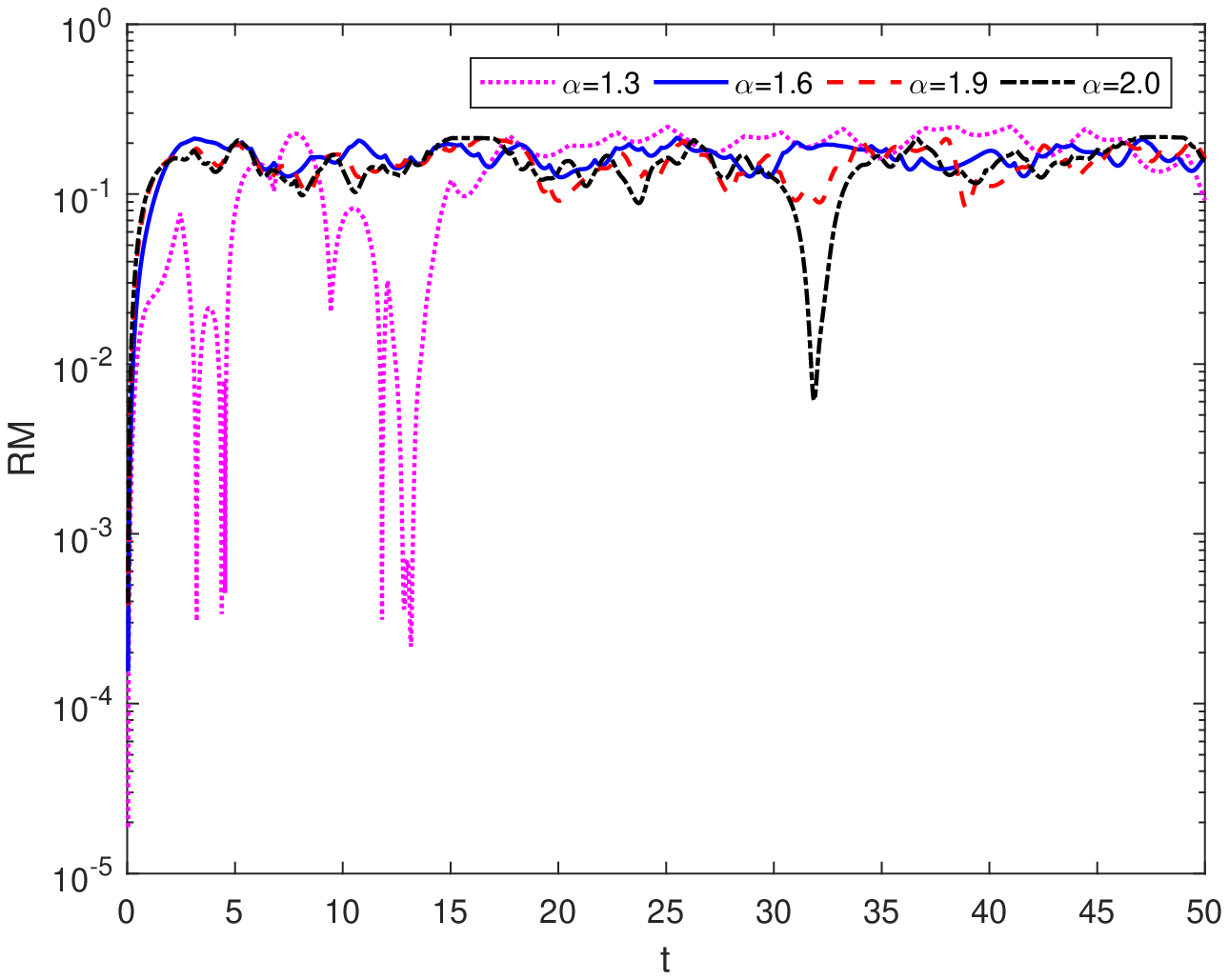}\\
{\footnotesize  \centerline {(b) Relative errors of mass}}
\end{minipage}
\caption{\small {\small{Relative errors of energy and mass with $h=0.5, \tau=0.02$ for different $\alpha$.} }}\label{fig522}
\end{figure}

Last but not least, we take the computation domain $\Omega=[-10,10]\times[-10,10]$ to investigate the relationship between the evolution of the soliton and the fractional order $\alpha$  for original system (\ref{NLS:eq:1.1}). The numerical results for different $\alpha$ with $N=256, \tau=0.001$ at $T=1$ are depicted in Figs. 5-8.
 Obviously, we can find that the $\alpha$ affect the shape of the wave function, the smaller the fractional order $\alpha$ is, the faster the wave function decays, and the steeper the shape of the wave function
  shape of the soliton changes dramatically. When the fractional order $\alpha$  tends to 2, the wave function converges to the classical NLS equation.

\begin{figure}[H]
\centering\begin{minipage}[t]{70mm}
\includegraphics[width=70mm]{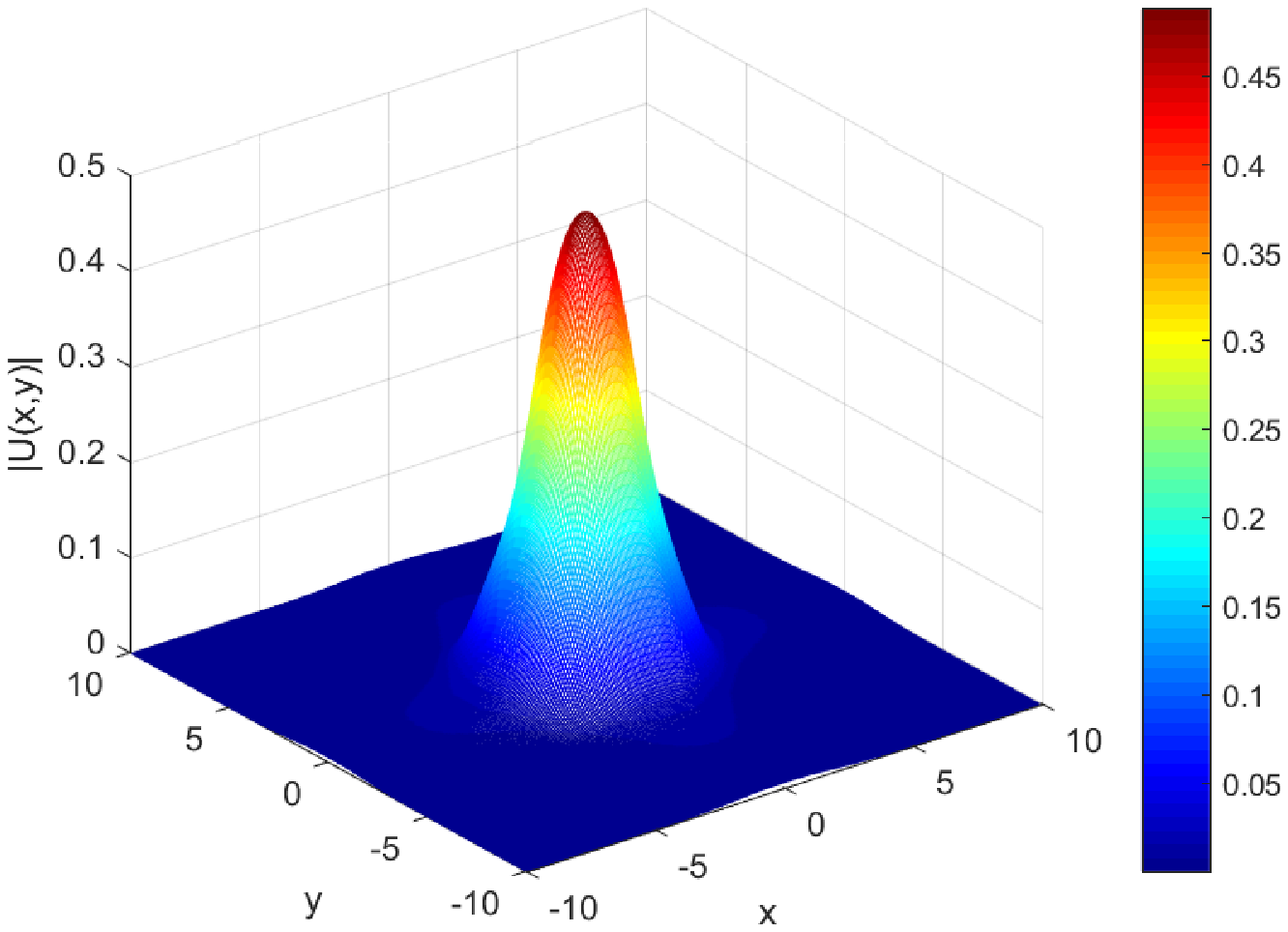}\\
\end{minipage}
\begin{minipage}[t]{70mm}
\includegraphics[width=70mm]{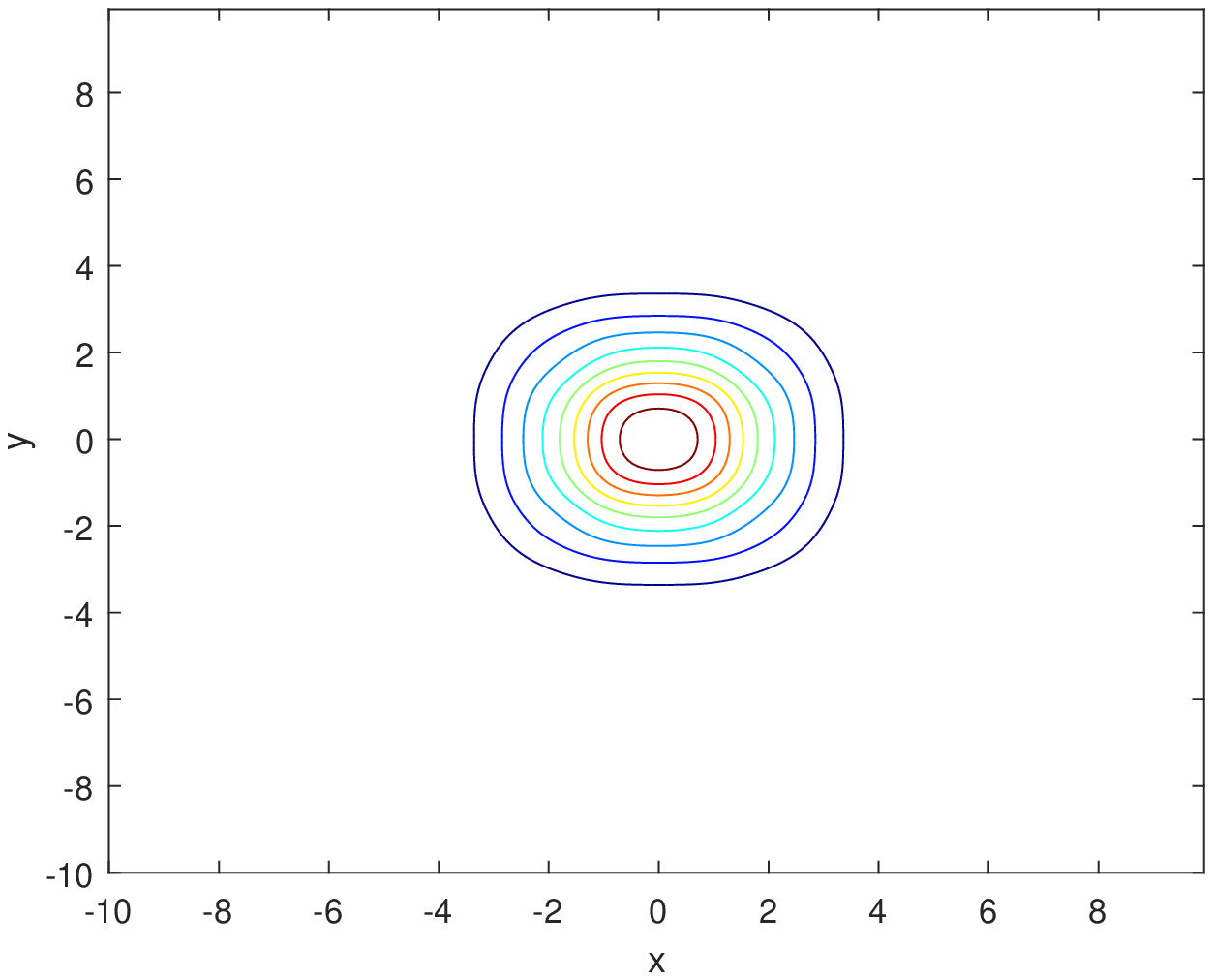}\\
\end{minipage}
\caption{\small {\small{Profile (left) and contour plot (right) of numerical solution for $\alpha=1.3$.} }}\label{fig522}
\end{figure}

\begin{figure}[H]
\centering\begin{minipage}[t]{70mm}
\includegraphics[width=70mm]{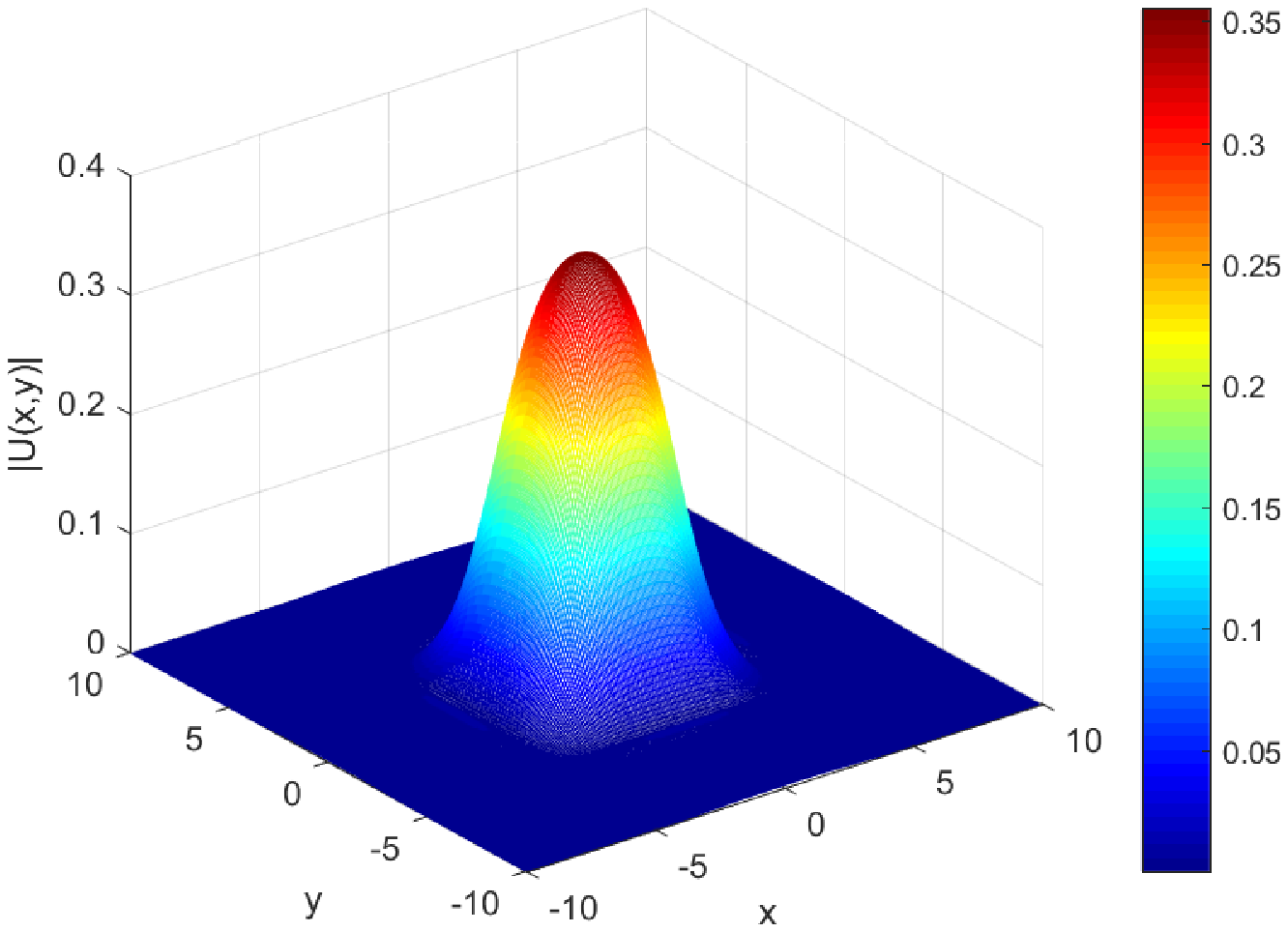}\\
\end{minipage}
\begin{minipage}[t]{70mm}
\includegraphics[width=70mm]{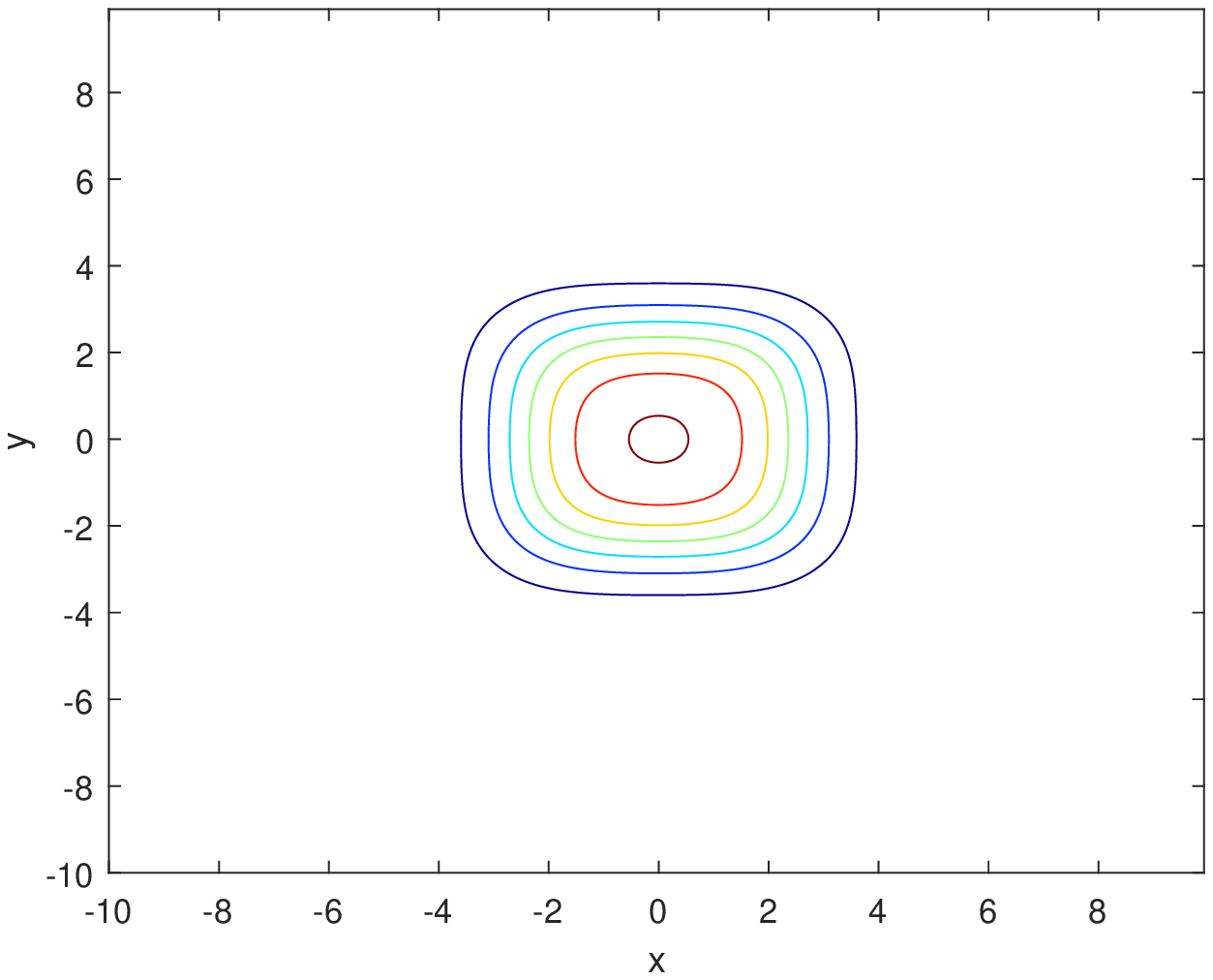}\\
\end{minipage}
\caption{\small {\small{Profile (left) and contour plot (right) of numerical solution for $\alpha=1.6$.} }}\label{fig522}
\end{figure}

\begin{figure}[H]
\centering\begin{minipage}[t]{70mm}
\includegraphics[width=70mm]{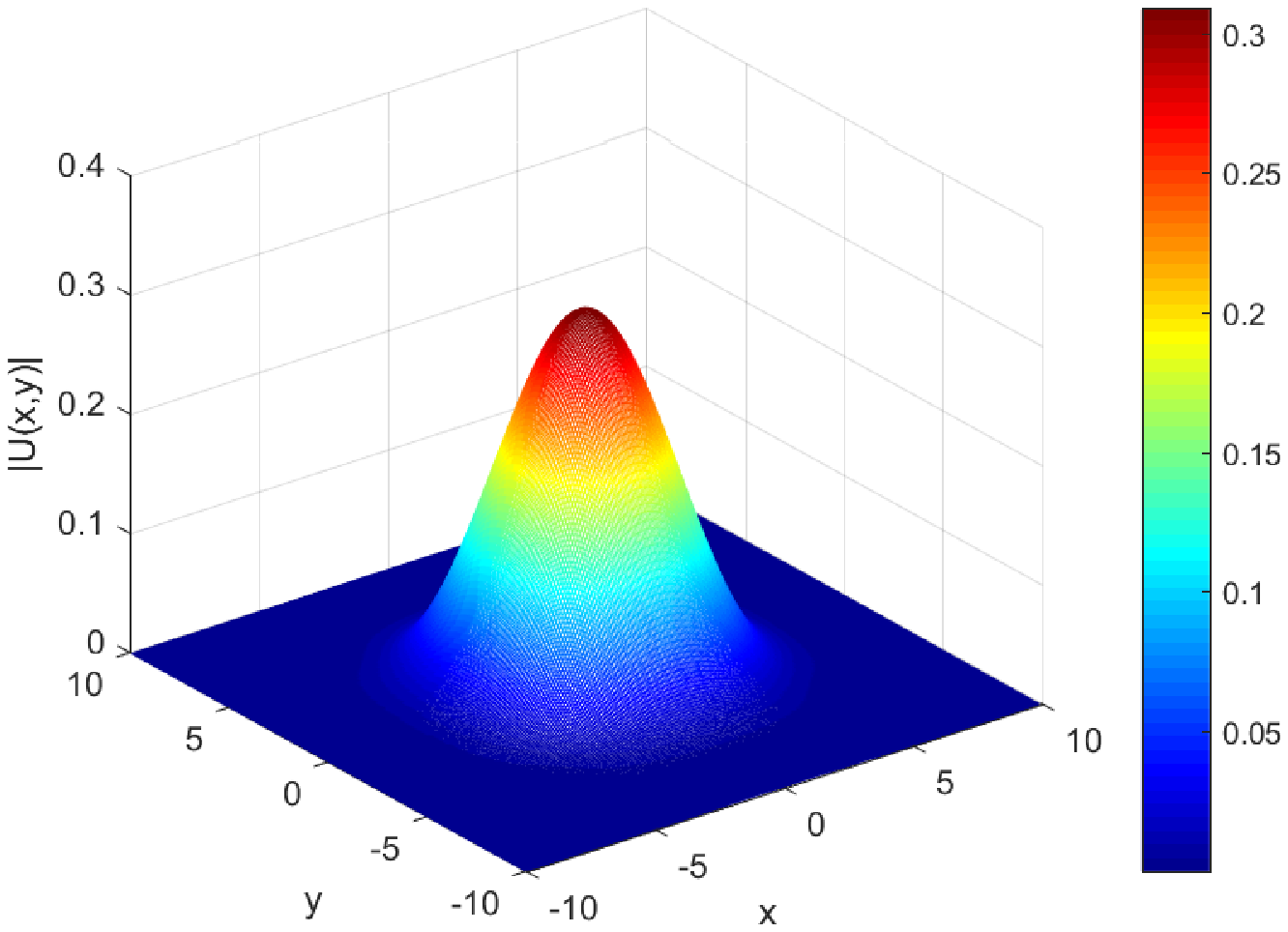}\\
\end{minipage}
\begin{minipage}[t]{70mm}
\includegraphics[width=70mm]{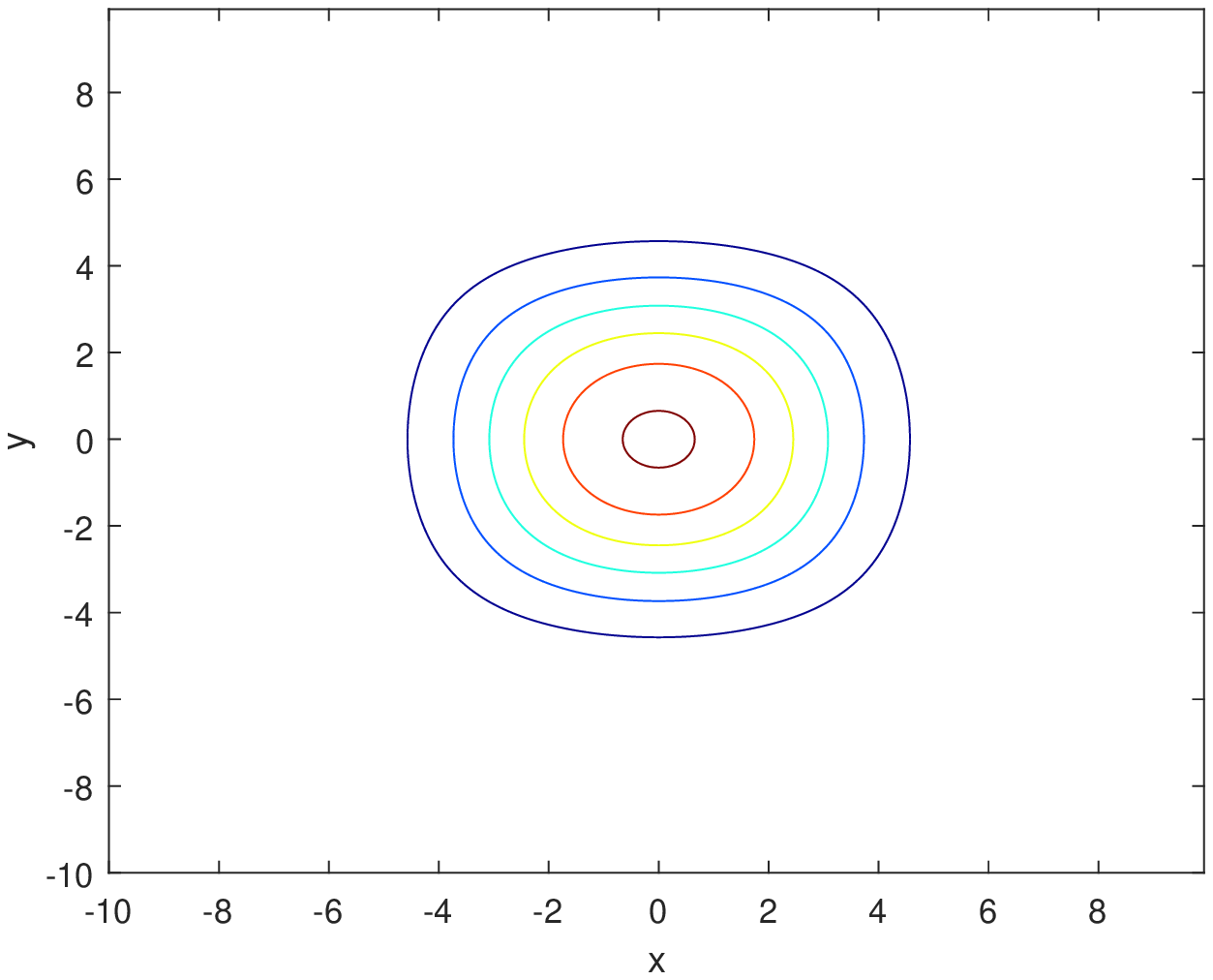}\\
\end{minipage}
\caption{\small {\small{Profile (left) and contour plot (right) of numerical solution for $\alpha=1.9$.} }}\label{fig522}
\end{figure}

\begin{figure}[H]
\centering\begin{minipage}[t]{70mm}
\includegraphics[width=70mm]{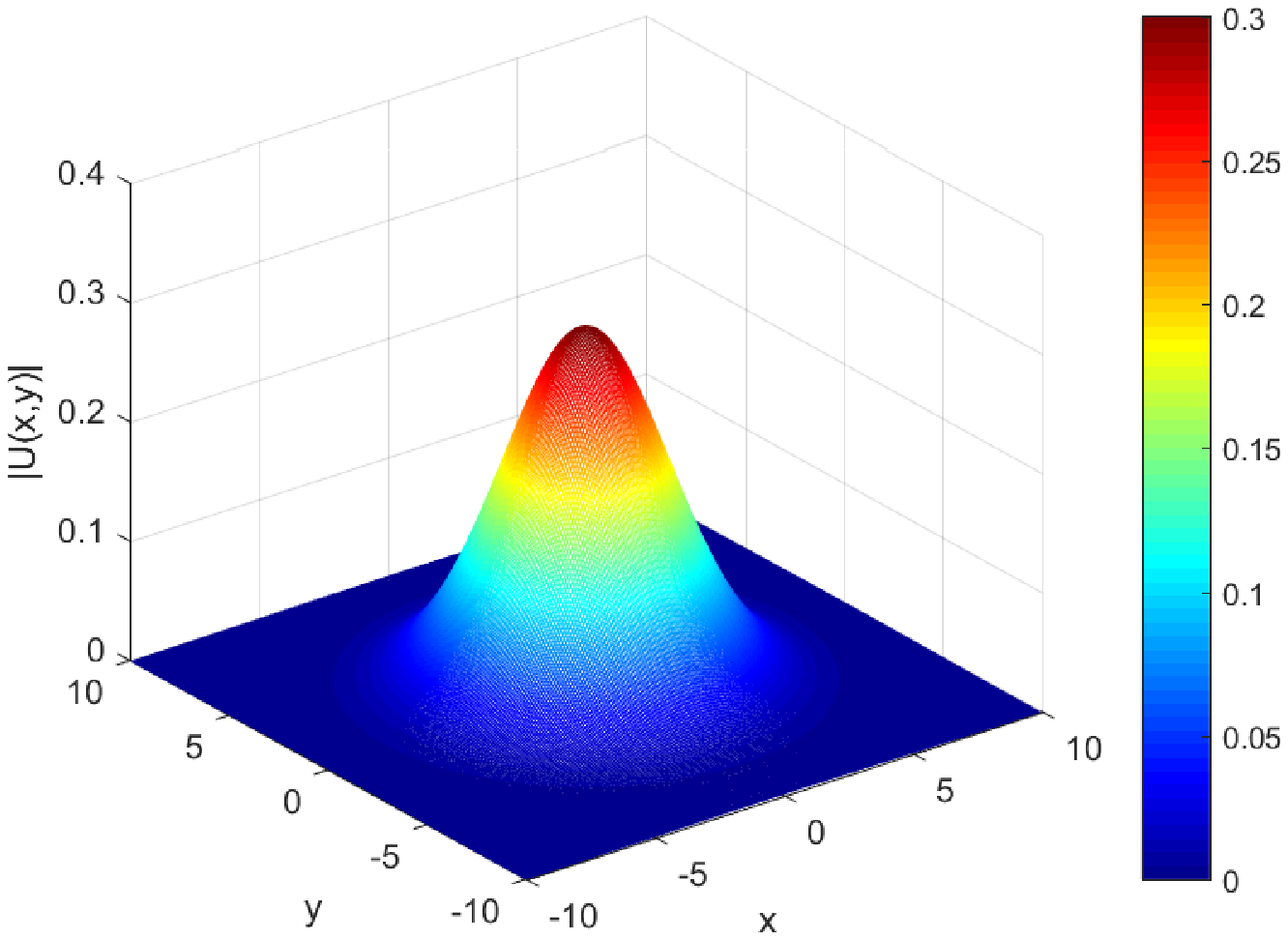}\\
\end{minipage}
\begin{minipage}[t]{70mm}
\includegraphics[width=70mm]{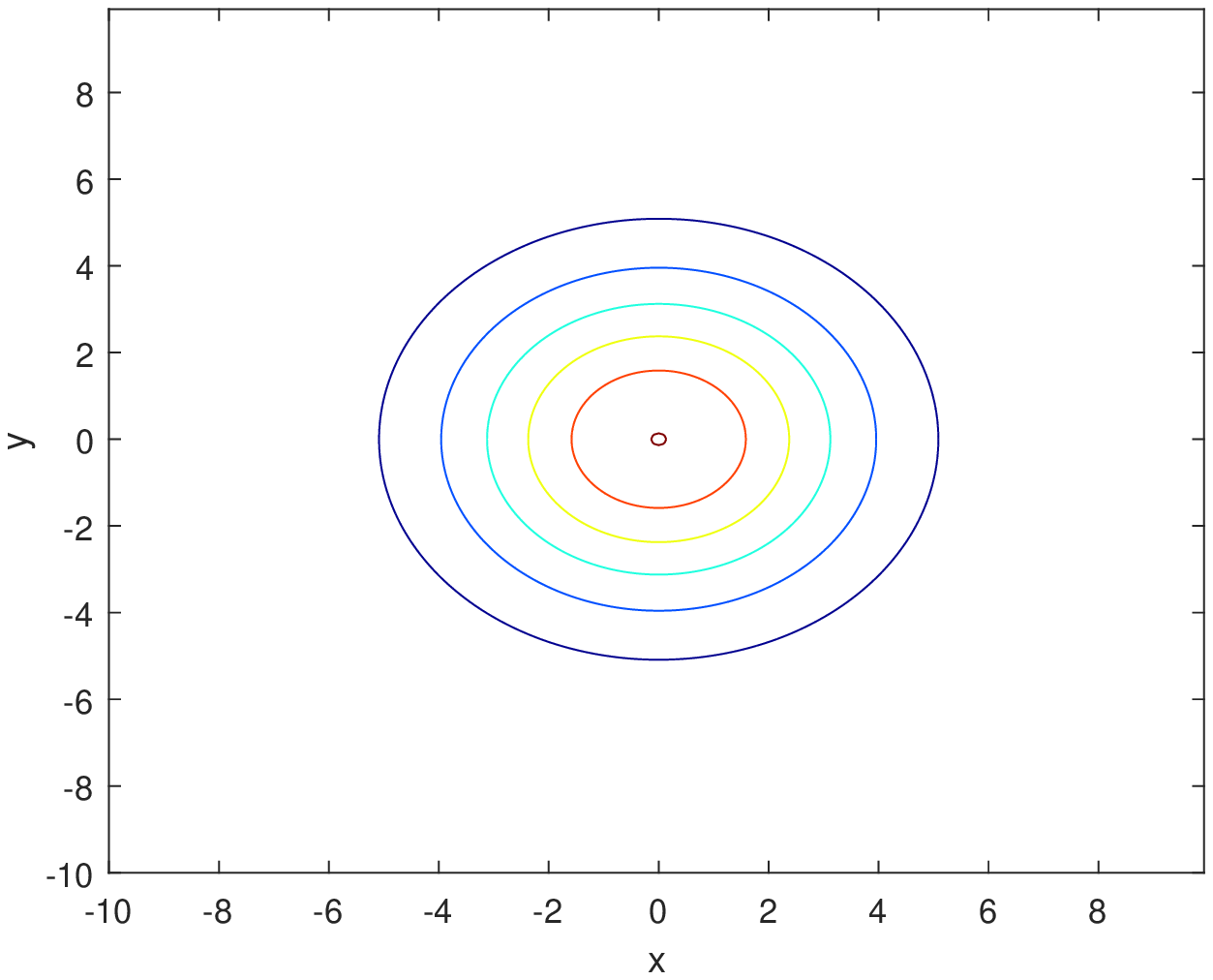}\\
\end{minipage}
\caption{\small {\small{Profile (left) and contour plot (right) of numerical solution for $\alpha=2$.} }}\label{fig522}
\end{figure}

\textbf{Example 4.3.} We consider the 2D fractional NLS equation (\ref{NLS:eq:1.1}) with potential functions. The initial condition is chosen to be
\begin{align}\label{NLS:eq:4.4}
u(x,y,0)=\frac{2}{\sqrt{\pi}}\text{exp}(-x^2-y^2).
\end{align}
In our computation, we take the computational domain $\Omega=[-5,5]\times[-5,5]$ with the parameters  $\gamma=\beta=1$,
 the harmonic potential $V_1$ and the optical lattice potential $V_2$ are given as
\begin{align*}
&V_1(x,y)=\frac{x^2+y^2}{2},\\
&V_2(x,y)=10(\text{sin}^2(\pi x)+\text{sin}^2(\pi y)).
\end{align*}

In Fig. 9, we show the relative energy errors for different potential functions $V$ when $h=0.5$,~$\tau=0.01$. Numerical results indicate that the FSAV scheme can preserve the discrete energy very well for the fractional NLS with a potential function. The contour plots of the numerical solutions of the fractional NLS with harmonic potential and optical lattice potential at different times are shown in Figs. 10-11, respectively.

\begin{figure}[H]
\centering\begin{minipage}[t]{70mm}
\includegraphics[width=70mm]{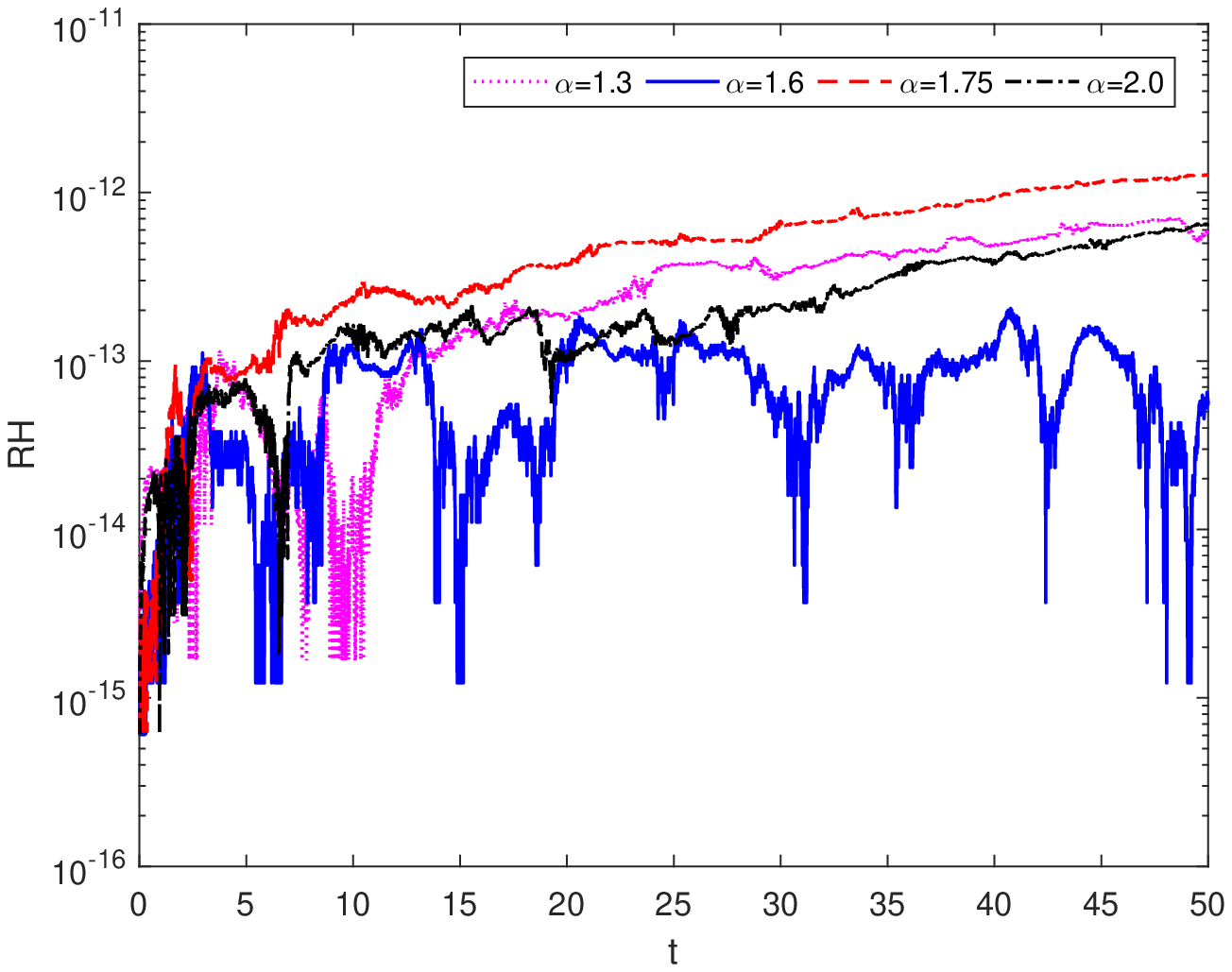}\\
{\footnotesize  \centerline {(a) \small{The relative energy errors for $V_1$.}}}
\end{minipage}
\begin{minipage}[t]{70mm}
\includegraphics[width=70mm]{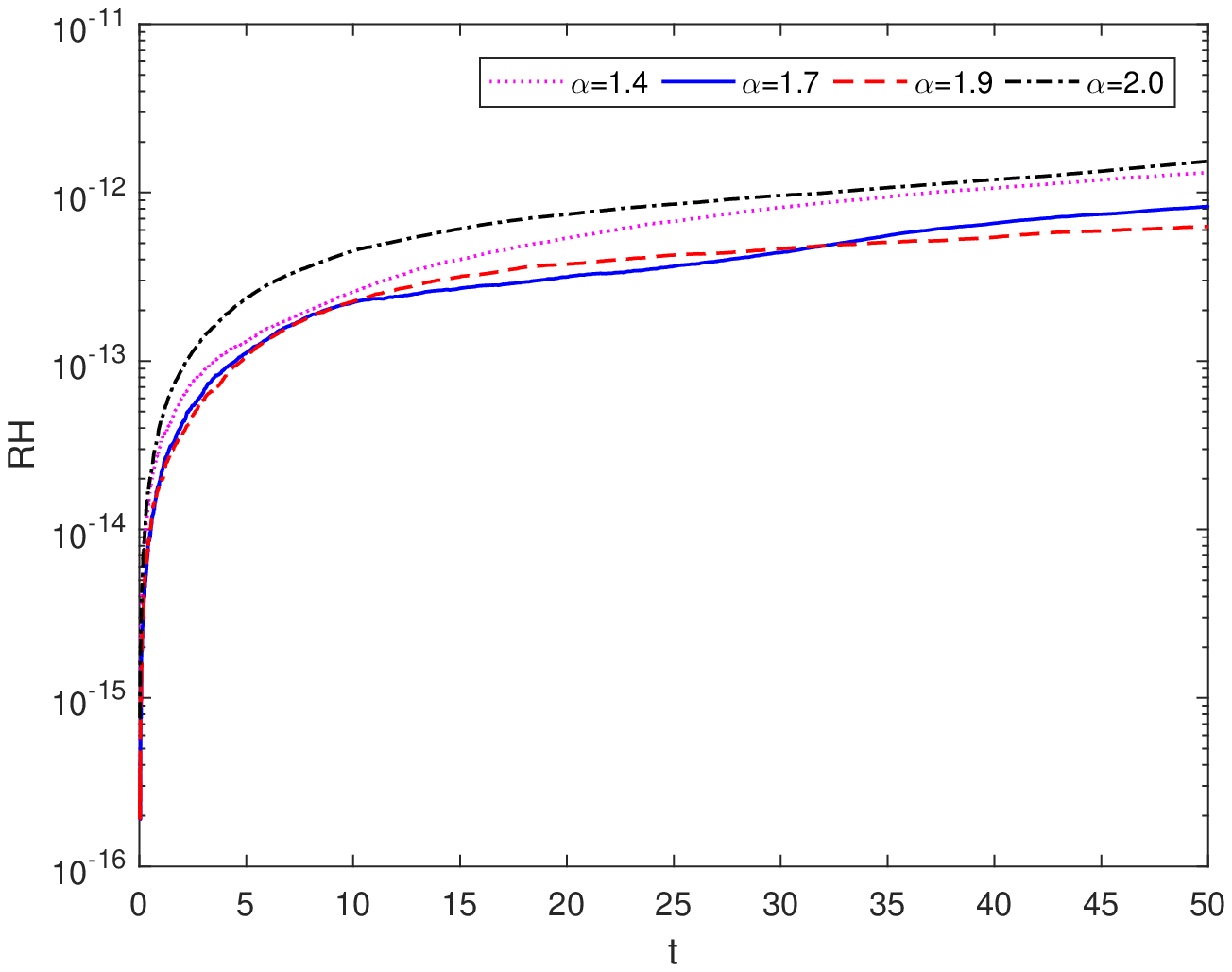}\\
{\footnotesize  \centerline {(b) \small{The relative energy errors for $V_2$ .}}}
\end{minipage}
\caption{\small {\small{The relative energy errors for different potential functions $V$ when $h=0.5$,~$\tau=0.01$.} }}\label{fig522}
\end{figure}

\begin{figure}[H]
\centering\begin{minipage}[t]{70mm}
\includegraphics[width=70mm]{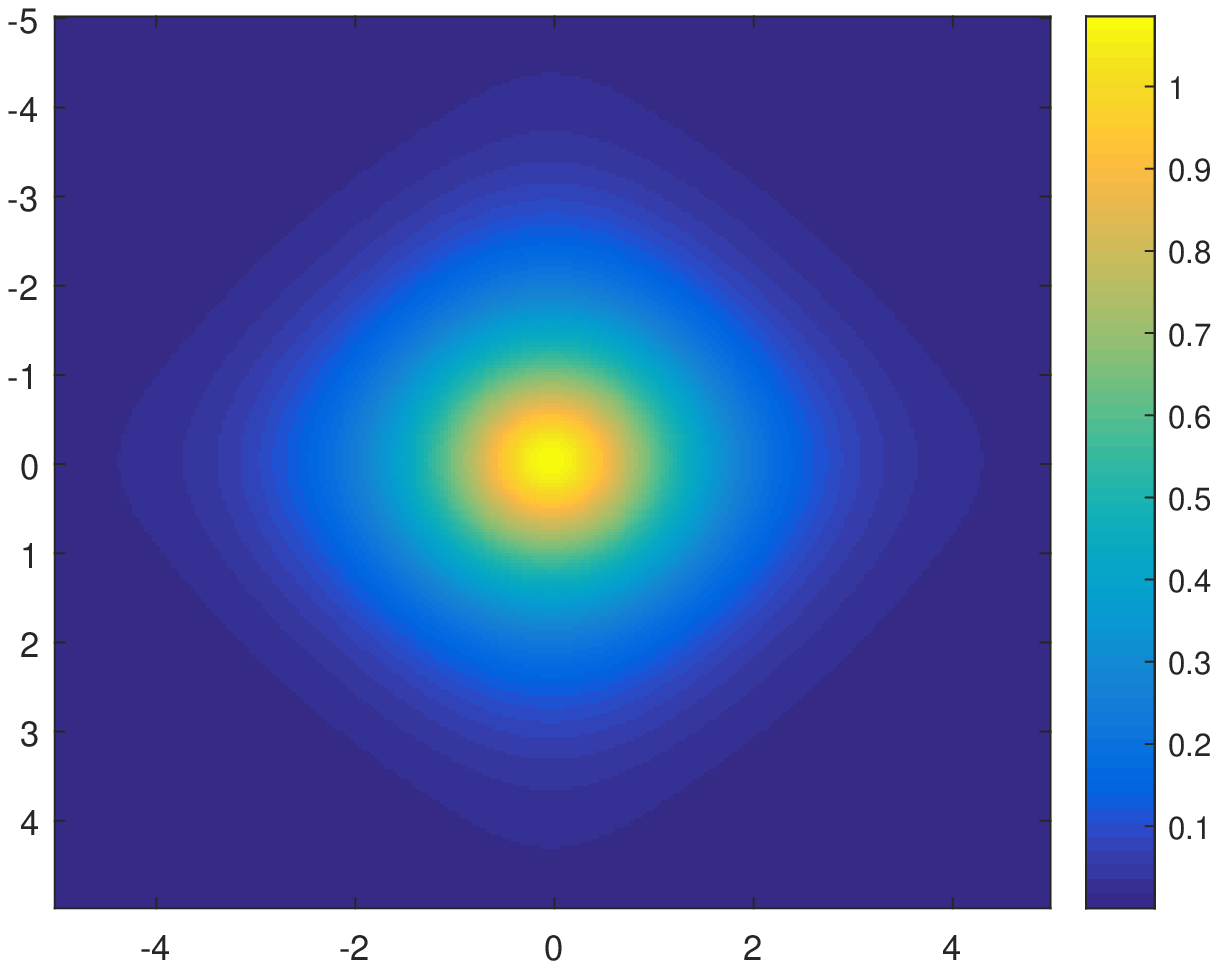}\\
{\footnotesize  \centerline {(a) $t=0.2$}}
\end{minipage}
\begin{minipage}[t]{70mm}
\includegraphics[width=70mm]{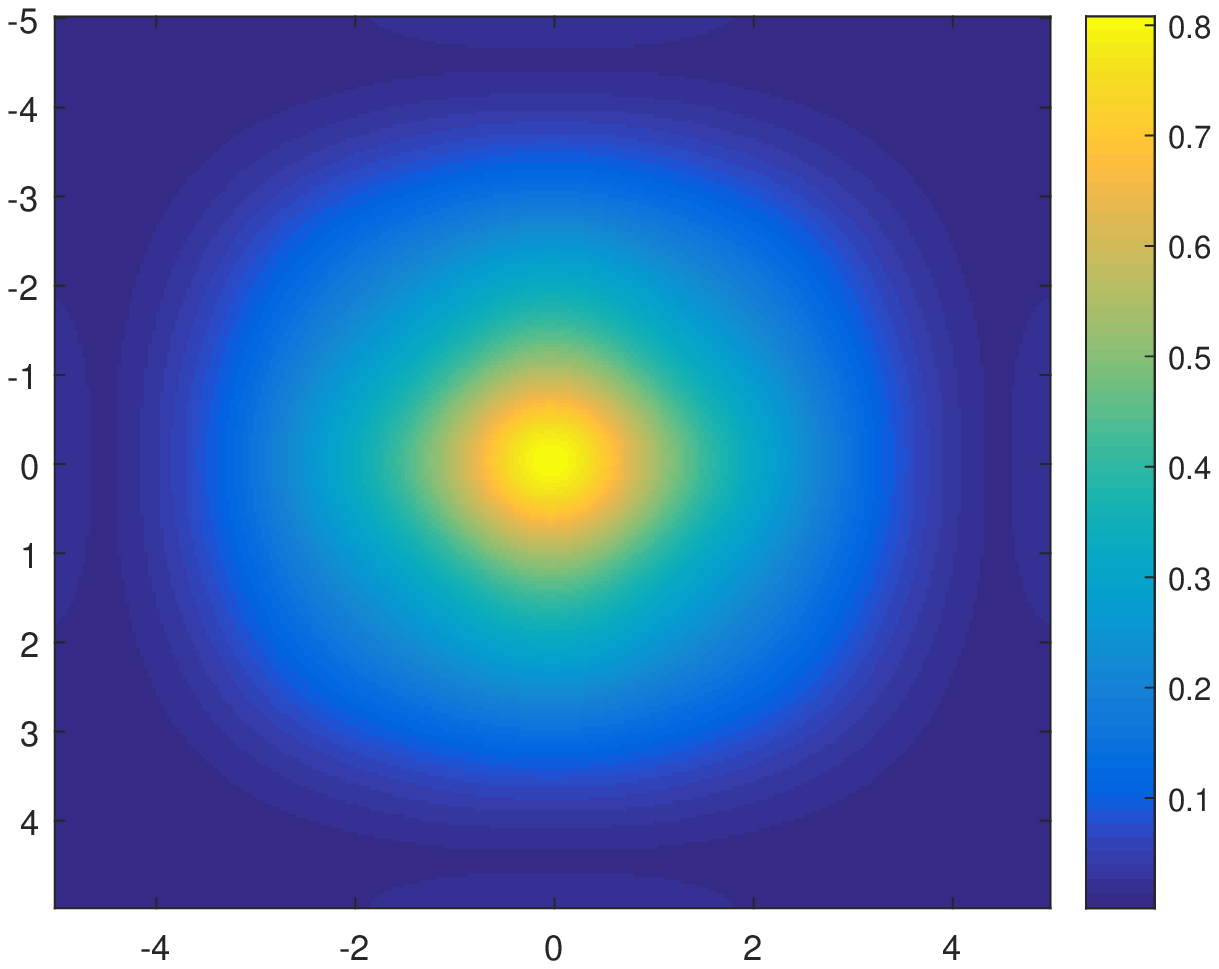}\\
{\footnotesize  \centerline {(b) $t=0.4$}}
\end{minipage}
\begin{minipage}[t]{70mm}
\includegraphics[width=70mm]{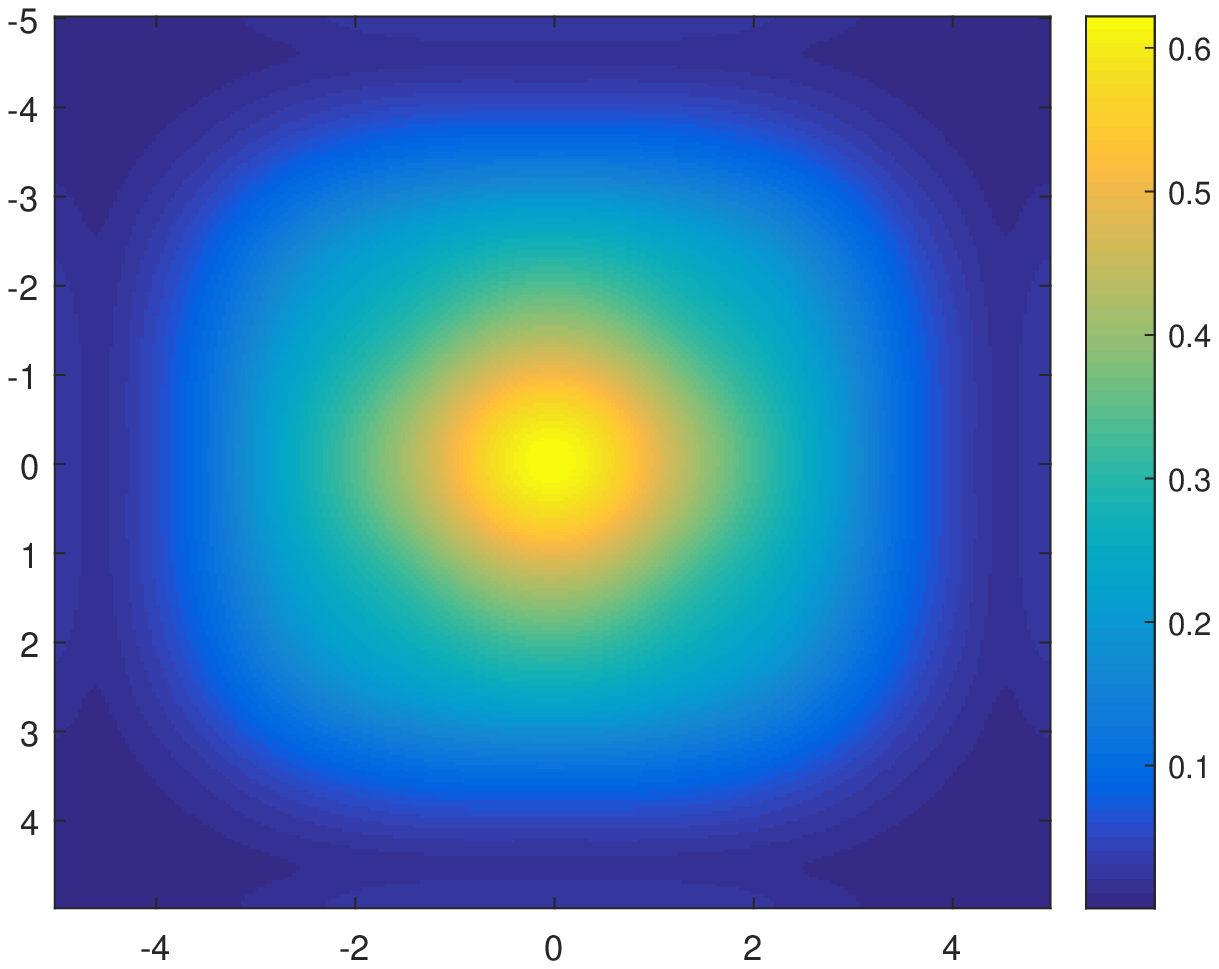}\\
{\footnotesize  \centerline {(c) $t=0.8$}}
\end{minipage}
\begin{minipage}[t]{70mm}
\includegraphics[width=70mm]{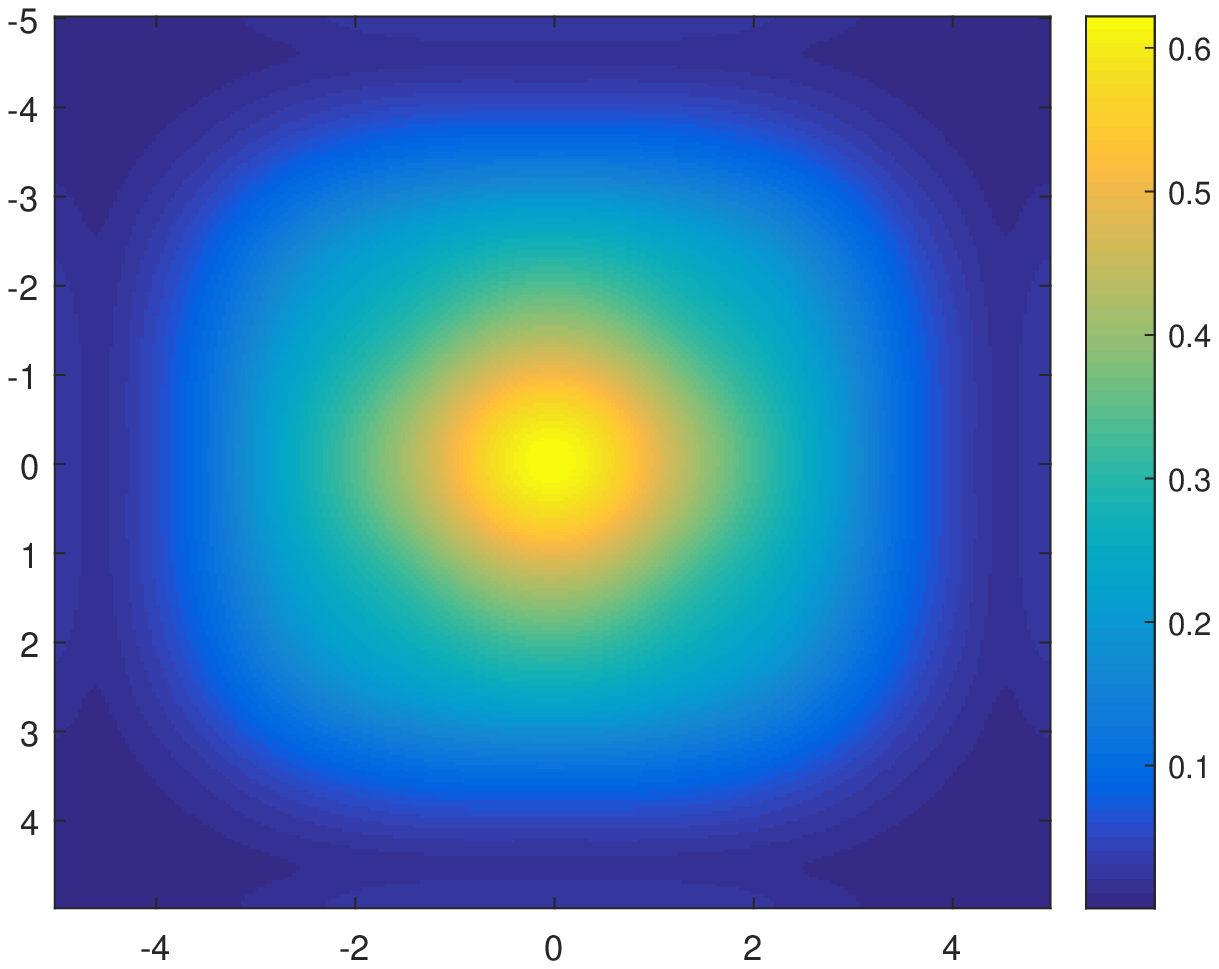}\\
{\footnotesize   \centerline {(d) $t=1.0$}}
\end{minipage}
\caption{\small {\small{Contour plots of the solution $|U(x,y)|$ for optical lattice potential $V_1$ at different times with $N=256, \tau=0.001$ for $\alpha=1.3$.} }}\label{fig522}
\end{figure}

\begin{figure}[H]
\centering\begin{minipage}[t]{70mm}
\includegraphics[width=70mm]{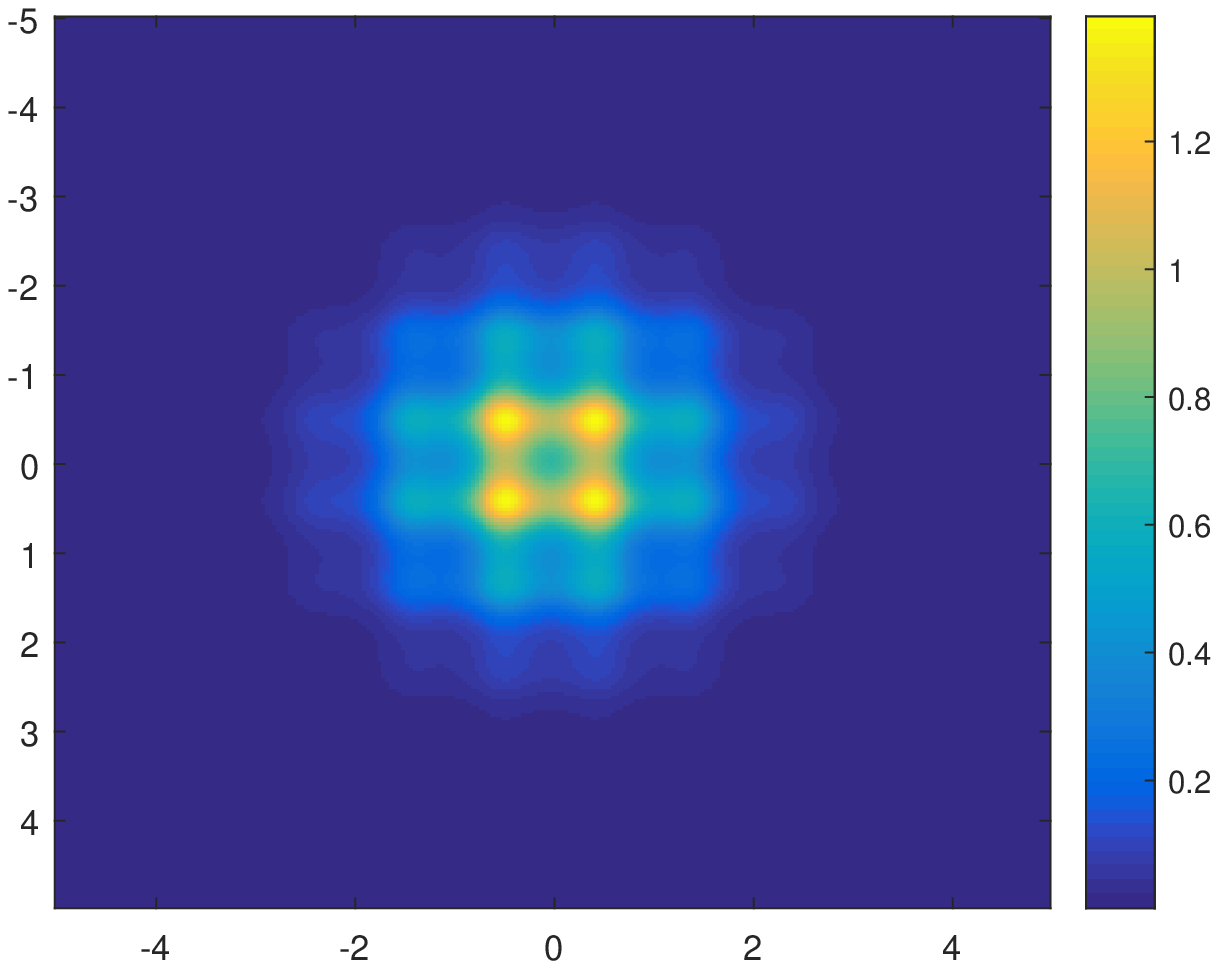}\\
{\footnotesize  \centerline {(a) $t=0.1$}}
\end{minipage}
\begin{minipage}[t]{70mm}
\includegraphics[width=70mm]{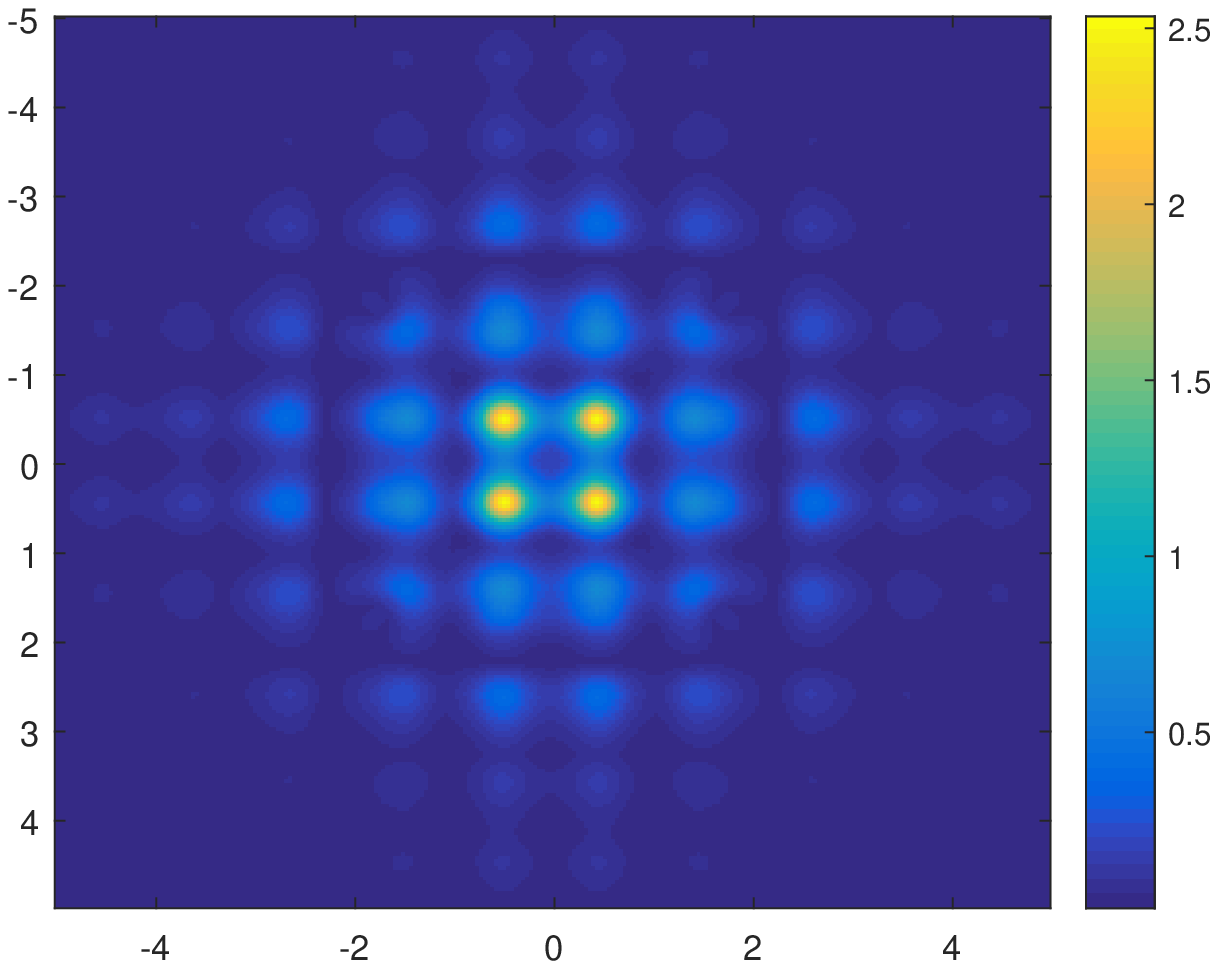}\\
{\footnotesize  \centerline {(b) $t=0.5$}}
\end{minipage}
\begin{minipage}[t]{70mm}
\includegraphics[width=70mm]{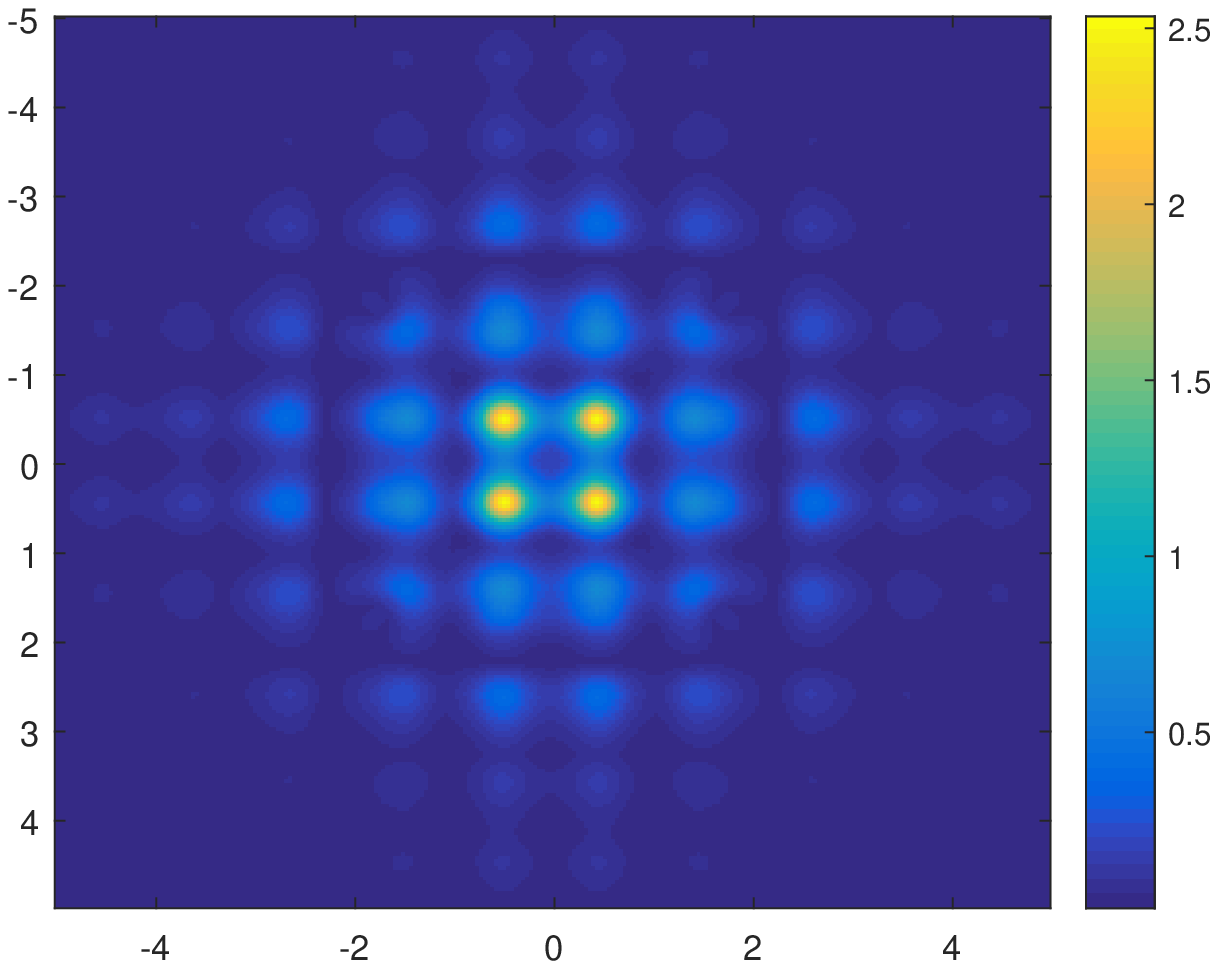}\\
{\footnotesize  \centerline {(c) $t=1.6$}}
\end{minipage}
\begin{minipage}[t]{70mm}
\includegraphics[width=70mm]{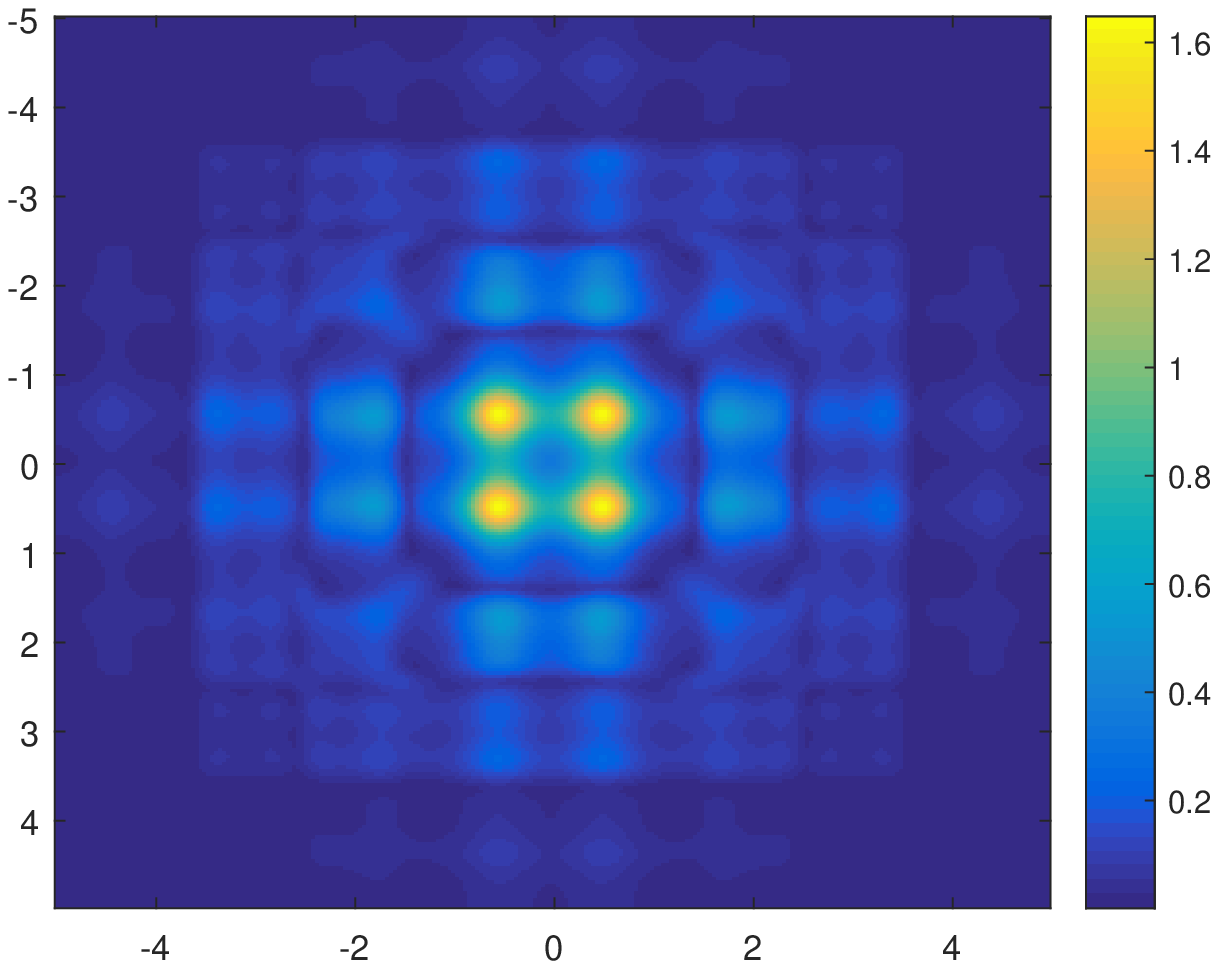}\\
{\footnotesize   \centerline {(d) $t=2.0$}}
\end{minipage}
\caption{\small {\small{Contour plots of the solution $|U(x,y)|$ for optical lattice potential $V_2$ at different times with $N=256, \tau=0.001$ for $\alpha=1.9$.} }}\label{fig522}
\end{figure}

\section{Conclusions}

%In this paper, a new linear Fourier pseudo-spectral scheme, being built upon the scalar auxiliary variable approach, is developed to solve the two-dimensional fractional nonlinear Schr\"{o}dinger equation. The presented scheme is linear, accurate and can conserve energy very well. Theoretical analysis and numerical results show that the new scheme is not only much more efficient and easy to implement, it also preserve energy conservation law. Specifically, the scheme here can be carried over to other fractional equations, such as the fractional sine-Gordon equation, the fractional Klein-Gordon-Schr\"{o}dinger equation, etc.. Future research could continue to explore analysing errors of the schemes which are based on scalar auxiliary variable approach approach.
%is not only much more efficient and easy to implement, it also
%constructing some liner conservative schemes, based on the scalar auxiliary variable approach, for solving fractional equations, such as the fractional Klein-Gordon-Schr\"{o}dinger equation, the fractional nonlinear wave equations, etc.

In this paper, taking the fractional nonlinear Schr\"{o}dinger equation as an example, we show that the scalar auxiliary variable approach can be extended to solve fractional differential equations. A linear implicit energy-preserving scheme is developed for the fractional nonlinear Schr\"{o}dinger equation. The energy conservation property and high efficiency of the proposed scheme are supported by theoretical analysis and numerical results. Following the similar process, the approach can be generalized to construct linear implicit energy-preserving scheme for solving other fractional differential equations, such as the fractional nonlinear wave equation, the fractional Klein-Gordon-Schr\"{o}dinger equation, etc. %Future research could continue to explore analysing errors of the numerical schemes that are based on the scalar auxiliary variable approach.
\section*{Acknowledgments}

This work is supported by the Postgraduate Research $\&$ Practice
Innovation Program of Jiangsu Province (Grant Nos. KYCX19\_0776), the National Natural Science Foundation of China (Grant No. 11771213, 61872422),
the National Key Research and Development Project of China (Grant No. 2016YFC0600310, 2018YFC0603500, 2018YFC1504205),
the Major Projects of Natural Sciences of University in Jiangsu Province of China (Grant No. 15KJA110002, 18KJA110003), and the Priority Academic Program Development of Jiangsu Higher Education Institutions.


\begin{thebibliography}{40}

\bibitem{p1}
F. Liu, and K. Burrage. \newblock Novel techniques in parameter estimation for fractional dynamical models arising from biological systems. \newblock {\em Comput. Math. Appl.} 62 (2011) 822-833.

\bibitem{p2}
F. Liu, C. Yang, and K. Burrage. \newblock Numerical method and analytical technique of the modified anomalous subdiffusion equation with a nonlinear source term. \newblock {\em J. Comput. Appl. Math.} 231 (2009) 160-176.

\bibitem{p3}
R. Metzler, and J. Klafter. \newblock The random walk's guide to anomalous diffusion: a fractional dynamics approach. \newblock {\em Phys. Rep.} 339 (2000) 1-77.

\bibitem{p4}
F. Liu, V. Anh, and I. Turner. \newblock Numerical solution of the space fractional Fokker-Planck equation. \newblock {\em J. Comput. Appl. Math.} 16 (2004) 209-219.

\bibitem{p5}
N. Laskin. \newblock Fractional quantum mechanics and L\'{e}vy path integrals. \newblock {\em Phys. Lett. A.} 268 (2000) 298-305.

\bibitem{p6}
N. Laskin. \newblock Fractional Schr\"{o}dinger equation. \newblock {\em Phys. Rev.} 66 (2002) 249-264.


\bibitem{p7}
B. Guo, Y. Han, and J. Xin. \newblock Existence of the global smooth solution to the period boundary value problem of fractional nonlinear Schr\"{o}dinger equation. \newblock {\em Appl. Math. Comput.} 204 (2008) 468-477.


\bibitem{p8}
S. Longhi. \newblock Fractional Schr\"{o}dinger equation in optics. \newblock {\em Opt. Lett.} 40 (2015) 1117-1120.


\bibitem{p13}
P. Constantin, and J. Wu. \newblock An extension problem related to the fractional Laplacian. \newblock {\em Commun. Partial Differ. Equ.} 26 (2009) 159-180.


\bibitem{p18}
D. Wang, A. Xiao, and W. Yang. \newblock A linearly implicit conservative difference scheme for the space fractional coupled nonlinear Schr\"{o}dinger equations. \newblock {\em J. Comput. Phys.} 272 (2014) 644-655.

\bibitem{p19}
P. Wang, and C. Huang. \newblock A conservative linearized difference scheme for the nonlinear fractional Schr\"{o}dinger equation. \newblock {\em Numer. Algorithms.} 69 (2015) 625-641.


\bibitem{p21}
M. Ran, and C. Zhang. \newblock A conservative difference scheme for solving the strongly coupled nonlinear fractional Schr\"{o}dinger equations. \newblock {\em Commun. Nonlinear Sci. Numer. Simul.} 41 (2016) 64-83.

\bibitem{p22}
X. Zhao, Z. Sun, and H. Peng. A fourth-order compact ADI scheme for 2D nonlinear space fractional Schr\"{o}dinger equation. \newblock {\em SIAM J. Sci. Comput.} 36 (2014) A2865-A2886.


\bibitem{p23}
S. Duo, and Y. Zhang. \newblock  Mass-conservative Fourier spectral methods for solving the fractional nonlinear Schr\"{o}dinger equation. \newblock {\em Comput. Math. Appl.} 71 (2016) 2257-2271.


\bibitem{p24}
H. Zhang, X. Jiang, W. Chu, and W. Fan. \newblock Crank-Nicolson Fourier spectral methods for the space fractional nonlinear Schr\"{o}dinger equation and its parameter estimation. \newblock {\em Int. J. Comput. Math.} DOI: 10.1080/00207160.2018.1434515.

\bibitem{p25}
Y. Wang, L. Mei, Q. Li, and L. Bu. \newblock Split-step spectral Galerkin method for the two-dimensional nonlinear space-fractional Schr\"{o}dinger equation. \newblock {\em Appl. Numer. Math.} 136 (2019) 257-278


\bibitem{p27}
M. Li, C. Huang, and P. Wang. Galerkin finite element method for nonlinear fractional Schr\"{o}dinger equations. \newblock {\em Numer. Algorithms.} 74 (2017) 499-525.

\bibitem{p29}
M. Li, X. Gu, C. Huang, and M. Fei. \newblock A fast linearized conservative finite element method for the strongly coupled nonlinear fractional Schr\"{o}dinger equations. \newblock {\em J. Comput. Phys.} 358 (2018) 256-282.

\bibitem{p30}
K. Feng, and M. Qin. \newblock Symplectic Geometric Algorithms for Hamiltonian Systems. \newblock {\em Springer Berlin Heidelberg.} 2010.

\bibitem{p31}
 E. Hairer,  C. Lubich, and G. Wanner. \newblock Geometric numerical integration: Structure-preserving algorithms for ordinary
differential equations. \newblock {\em 2nd edition, Springer-Verlag, Berlin}, 2006.

\bibitem{p32}
B. Leimkuhler, and S. Reich. Simulating Hamiltonian Dynamics. \newblock {\em Cambridge University Press, Cambridge,} 2004.

\bibitem{p33}
F. Zhang, V. P\'{e}rez-Garc\'{\i}a, and L. V\'{a}zquez. \newblock Numerical simulation of nonlinear Schr\"{o}dinger systems: A new conservative scheme. \newblock {\em Appl. Math. Comput.} 71 (1995) 165-177.

\bibitem{p34}
P. Wang, and C. Huang. \newblock Structure-preserving numerical methods for the fractional Schr\"{o}dinger equation. \newblock {\em Appl. Numer. Math.} 129 (2018) 137-158.

\bibitem{p35}
P. Wang, and C. Huang. \newblock Split-step alternating direction implicit difference scheme for the fractional Schr\"{o}dinger equation in two dimensions. \newblock {\em Comput. Math. Appl.} 71 (2016) 1114-1128.

\bibitem{p36}
X. Yang, J. Zhao, and Q. Wang. \newblock Linear, first and second-order, unconditionally energy stable numerical
schemes for the phase field model of homopolymer blends. \newblock {\em J. Comput. Phys.} 327 (2016) 294-316.


\bibitem{p37}
Y. Gong, J. Zhao, X. Yang, and Q. Wang. \newblock Fully discrete second-order linear schemes for hydrodynamic phase field models of binary viscous fluid flows with variable densities. \newblock {\em SIAM J. Sci. Comput.} 40 (2018) B138-B167.

\bibitem{p38}
J. Shen, J, Xu, and J, Yang. \newblock The scalar auxiliary variable (SAV) approach for gradient flows. \newblock {\em J. Comput. Phys.} 353 (2018) 407-416.

\bibitem{p39}
J. Shen, J, Xu, and J, Yang. \newblock A new class of efficient and robust energy stable schemes for gradientows. \newblock{\em arXiv:1710.01331,} 2017.

\bibitem{p40}
J. Mac\'{\i}as-D\'{\i}az. \newblock A structure-preserving method for a class of nonlinear dissipative wave equations with Riesz space-fractional derivatives. \newblock {\em J. Comput. Phys.} 351 (2017) 40-58.

\bibitem{p41}
J. Mac\'{\i}as-D\'{\i}az, and E. Jorge. \newblock A numerically efficient dissipation-preserving implicit method for a nonlinear multidimensional fractional wave equation. \newblock {\em J. Sci. Comput.} 77 (2018) 1-26.

\bibitem{p42}
 J. Wang, and A, Xiao. \newblock An efficient conservative difference scheme for fractional Klein-Gordon-Schr\"{o}dinger equations. \newblock {\em Appl. Math. Comput.} 320 (2018) 691-709.

\bibitem{p43}
S. Secchi, and M. Squassina. \newblock Soliton dynamics for fractional Schr\"{o}dinger equations. \newblock {\em Appl. Anal.} 93 (2014) 1702-1729.

\bibitem{p44}
M. Cheng. \newblock Bound state for the fractional Schr\"{o}dinger equation with unbounded potential. \newblock {\em J. Math. Phys.} 53 (2012) 043507.

\bibitem{p45}
J. Shen and T. Tang. \newblock Spectral and High-Order Methods with Applications. \newblock {\em Science Press, Beijing.} 2006.

\bibitem{p46}
Q. Cheng, J. Shen and X. Yang. \newblock Highly efficient and accurate numerical schemes for the epitaxial thin film growth models by using the SAV approach. \newblock {\em J. Sci. Comput.} DOI: 10.1007/s10915-018-0832-5.


\end{thebibliography}
\end{document}